\documentclass[11pt,reqno]{article}

\usepackage{amsmath,amsthm,amsfonts,amssymb,bm,enumerate,color,tikz,hyperref}
\usetikzlibrary{matrix,arrows,positioning,calc}
\usepackage{tikzpagenodes}
\usepackage[margin=1.25in]{geometry}
\usepackage{etex}
\usepackage{cite}

\usepackage{pifont}
\usepackage{makecell} 
\usepackage{multirow}
\usepackage{latexsym}
\usepackage{mathtools}

\usepackage{ragged2e}
\usepackage{graphics,graphicx,psfrag}
\usepackage{amscd}
\usepackage{epsfig}
\usepackage[all]{xy}
\usepackage{float}
\usepackage{enumitem}
\usetikzlibrary{decorations.pathmorphing}
\usepackage{setspace}

\newtheoremstyle{break}
  {0.5cm}%
  {0.5cm}%
  {\itshape}%
  {}%
  {\bfseries}%
  {\vspace{3pt}}%
  {\newline}%
  {}%
\theoremstyle{break}

\newtheorem{definition}{Definition}[section]
\newtheorem{theorem}[definition]{Theorem}
\newtheorem{lemma}[definition]{Lemma}
\newtheorem{corollary}[definition]{Corollary}
\newtheorem{proposition}[definition]{Proposition}

\newtheorem{example}[definition]{Example}

\newenvironment{myproof}[1][\proofname]{%
  \begin{proof}[\textbf{\textit{Proof}}]$ $\par\vspace{-2pt}
}{%
  \end{proof}
}

\numberwithin{equation}{section}

\newcommand{\Mor}{\mathsf{Mor}}
\newcommand{\Hom}{\mathsf{Hom}}

\newcommand{\too}{\longrightarrow}

\newcommand{\A}{\mathcal{A}}
\newcommand{\B}{\mathcal{B}}

\newcommand{\C}{\mathcal{C}}
\newcommand{\D}{\mathcal{D}}
\newcommand{\E}{\mathcal{E}}
\newcommand{\F}{\mathcal{F}}

\newcommand{\Mod}{\mathsf{Mod}}
\newcommand{\Ker}{\mathrm{Ker}}

\newcommand{\Coker}{{\rm CoKer}}
\newcommand{\op}{{}^{\rm op}}
\newcommand{\Set}{\mathsf{Set}}
\newcommand{\Cat}{\mathsf{Cat}}
\newcommand{\Ab}{\mathsf{Ab}}
\newcommand{\Fun}{\mathsf{Fun}}
\newcommand{\Lan}{{\rm Lan}}

\newcommand{\uno}{\mathsf{(1)}}
\newcommand{\dos}{\mathsf{(2)}}
\newcommand{\tres}{\mathsf{(3)}}

\newcommand{\iroman}{\mathsf{(i)}}
\newcommand{\iiroman}{\mathsf{(ii)}}

\hypersetup{
    bookmarks=true,         
    unicode=false,          
    pdftoolbar=true,        
    pdfmenubar=true,        
    pdffitwindow=false,     
    pdfstartview={FitH},    
    pdftitle={My title},    
    pdfauthor={Author},     
    pdfsubject={Subject},   
    pdfcreator={Creator},   
    pdfproducer={Producer}, 
    pdfkeywords={keyword1} {key2} {key3}, 
    pdfnewwindow=true,      
    colorlinks=true,       
    linkcolor=black,          
    citecolor=magenta,        
    filecolor=cyan,      
    urlcolor=violet           
}

\usepackage{fancyhdr}

\fancypagestyle{plain}{%
\fancyhf{}
\fancyhead[L]{M. A. P\'erez}\fancyhead[C]{}\fancyhead[R]{Coends and the tensor product of $\C$-modules}
\fancyfoot[L]{}\fancyfoot[C]{}\fancyfoot[R]{Page \thepage}
\pagestyle{fancy}
\thispagestyle{plain}
}

\usepackage{tocloft}

\usepackage{titlesec}

\titleformat{\subsubsection}[runin]
  {\normalfont\bfseries}{\thesubsection}{1em}{}

\begin{document}

{\huge Coends and the tensor product of $\C$-modules} 

\vspace{0cm}

\vskip3mm \noindent Marco A. P\'erez \\ 
Instituto de Matem\'aticas \\ 
Universidad Nacional Aut\'onoma de M\'exico \\ 
Circuito Exterior, Ciudad Universitaria \\
C.P. 04510, Mexico City, MEXICO \\ 
{\tt maperez@im.unam.mx}

\vfill

{\small The author is supported by a DGAPA postdoctoral fellowship from the Universidad Nacional Aut\'onoma de M\'exico}

\

\hfill \today

\thispagestyle{empty}

\newpage


~\

\vspace{6.5cm}

\thispagestyle{empty}

\begin{abstract} 
\noindent We give an introduction to the concept of Kan extensions, and study its relation with the notions of coend and adjoint functors. We state and prove in detail a well known formula to compute Kan extensions by using coends: a certain colimit related to the concept of copower. Finally, we study the tensor product of functors, and its relation with Kan extensions, in order to represent the tensor product of $\C$-modules as a particular case. 
\end{abstract}

\newpage


\tableofcontents

\listoffigures

\newpage


\section*{Introduction}\label{sec:intro}

\addcontentsline{toc}{section}{Introduction}

The notion of tensor product is ubiquitous in mathematics, but also in other fields such as Physics. The reader probably had her/his first exposure to this concept in her/his first courses of linear algebra, when studied tensor product of vector spaces. People who specialized in algebra in their undergraduate mathematics program took one step forward to the generalization of tensor product of vector spaces, by studying the tensor product of modules. In the branch known as homological algebra, it is usual to study at least two different kinds of tensor products of chain complexes. This point represented for the author a first contact with the notion of \emph{monoidal category}, where tensor products are studied in the abstract framework provided by category theory. The examples of tensor products mentioned so far satisfy certain \emph{universal property}, and so it is not surprising that tensor products appear as a construction widely studied in category theory. Categorical settings give rise to applications in other branches of mathematics, such as in representation theory of Artin algebras (specifically, in tilting theory). In this particular field, the notion of tensor product of $\C$-modules (with $\C$ a \emph{skeletally small} and preadditive category) plays an important role. The main purpose of these notes is to give to this tensor product an interpretation as a tensor product in a certain category of functors, but for such a task, it comes handy to be familiar with the concept of Kan extensions.

Besides generalizing the concept of colimit, Kan extensions are universal constructions related to the notions of adjoint functors and coends. In \emph{Section~\ref{sec:def}}, we recall the definition of Kan extensions and point out its universality by studying their relation with colimits and adjoint functors. \emph{Section~\ref{sec:formula}} is devoted to study the interplay between Kan extensions and coends. We will see that under certain conditions, it is possible to compute Kan extensions via a formula involving the coend of a certain bifunctor (See \emph{Theorem~\ref{theo:MacLane}}). We will give a detailed proof of this formula, which only demands a basic knowledge of universal constructions, namely, coproducts and coequalizers. Finally, in \emph{Section~\ref{sec:producto}} we will recall the concept of tensor product of functors as a particular coend, and will show how to relate it to the tensor product of $\C$-modules, widely used in tilting theory. In order to be able to establish this relation, we need to keep in mind that given a $\C\op$-module $F \colon \C \too \Ab$ and a $\C$-module $G \colon \C\op \too \Ab$ their tensor product (of functors) $F \otimes_{\C\op} G$ exists in the case where $\C$ is small. However, we will see that this fact can be extended to skeletally small categories, by showing a slight generalization of \emph{Theorem~\ref{theo:MacLane}} (See \emph{Corollary~\ref{theo:MacLaneExtension}}).


\subsection*{Notation and preliminaries}

Throughout the present notes, we assume the reader is familiar with the following concepts from category theory: universal constructions in categories, functors, natural transformations, representable functors, and adjoint functors; as well as, kernels, cokernels, and exact sequences in abelian categories. A good introduction to these topics is the book \emph{\cite{Leinster}} by T. Leinster. We also assume that the reader is familiar with the notion of \emph{duality}. All concepts and results presented in these notes have their corresponding dual versions, which will be omitted for the sake of simplicity. 

Given a category $\C$, denote by $\Mor(\C)$ the class of morphisms in $\C$. By $X \in \C$ we will mean that $X$ is an object of $\C$. Given $X, Y \in \C$, we denote by $\Hom_{\C}(X,Y)$ the collection of morphisms from $X$ to $Y$. If two objects $X$ and $Y$ are isomorphic, we will use the notation $X \simeq Y$. A category $\C'$ is a \textbf{subcategory} of $\C$ if every object of $\C'$ is an object of $\C$, and if every morphism of $\C'$ is a morphism of $\C$. If in addition, for every $X \in \C$ there exists $X' \in \C'$ such that $X \simeq X'$, $\C'$ is called \textbf{dense} in $\C$.

For the purpose of these notes, $\Hom_{\C}(X,Y)$ will always be a set, for every $X, Y \in \C$. In other words, categories considered from now on will be \textbf{locally small}. If the collection of objects of $\C$ is also a set, then $\C$ is what is known as a \textbf{small category}. Between small and locally small categories, roughly speaking, we have the notion of \textbf{skeletally small} category, that is, a category $\C$ containing a subcategory $\C'$ which is dense in $\C$. 

We will use short arrows $\to$ to represent morphisms in a category. Thus, $f \colon X \to Y$ denotes a morphism in a category, say $\C$, from the object $X$ to the object $Y$. Long arrows $\too$ will be used to represent functors between categories; for example, $F \colon \C \too \D$ denotes a functor from the category $\C$ to the category $\D$. Finally, double arrows $\Rightarrow$ will represent natural transformations between functors: given two functors $F, G \colon \C \too \D$, then $\alpha \colon F \Rightarrow G$ denotes a natural transformation from $F$ to $G$. In the case $\alpha$ is a natural isomorphism, we will write $F \cong G$. 

Given a natural transformation $\alpha \colon F \Rightarrow G$ between functors $F, G \colon \C \too \D$, and two functors $K \colon \A \too \C$ and $H \colon \D \too \B$, we define:
\begin{itemize}[itemsep=2pt,topsep=0pt]
\item[$\uno$] $\alpha_K \colon F \circ K \Rightarrow G \circ K$ as the natural transformation given by the family of morphisms
\[
(\alpha_K)_A := \alpha_{K(A)}, \mbox{ for every }A \in \A.
\] 

\item[$\dos$] $H(\alpha) \colon H \circ F \Rightarrow H \circ G$ as the natural transformation given by the family of morphisms 
\[
(H(\alpha))_C := H(\alpha_C), \mbox{ for every } C \in \C.
\] 
\end{itemize}
There are more operations between natural transformations, known as \emph{horizontal} and \emph{vertical} compositions, which give a 2-category structure to the collection of small categories, functors, and natural transformations. For the purpose of the present notes, we only need to recall how the vertical composition is defined. The remaining operations and axioms defining a 2-category will be omitted, but we refer the interested reader to \emph{\cite{HigherLeinster}}, also by T. Leinster. 
\begin{itemize}[itemsep=2pt,topsep=0pt]
\item[$\tres$] Let $F, G, H \colon \C \too \D$ be functors between categories, and $\alpha \colon F \Rightarrow G$ and $\beta \colon G \Rightarrow H$ be natural transformations. The (vertical) composition $\beta \circ \alpha \colon F \Rightarrow G$ is the natural transformation formed by the following family of morphisms:
\[
(\beta \circ \alpha)_C := \beta_C \circ \alpha_C, \mbox{ for every } C \in \C.
\]
\end{itemize}
In certain occasions, we will work with some \emph{concrete} categories, namely: the category $\Set$ of sets and functions; the categories $\Mod(R)$ and $\Mod(R\op)$ of left and right $R$-modules, respectively, and $R$-homomorphisms; and the category $\Ab$ of abelian groups (that is, $\mathbb{Z}$-modules) and homomorphisms. Given two categories $\C$ and $\D$, we denote by $\Fun(\C,\D)$ the category whose objects are the functors from $\C$ to $\D$, and whose morphisms are the natural transformations between such functors.


\section{Kan extensions}\label{sec:def}

We begin this section recalling the concept of Kan extensions. For simplicity, and due to the interests of these notes, we only work with left Kan extensions, and omit their dual analogous of right Kan extensions. 

Suppose we are given the following diagram of functors
\begin{equation}\label{fig:extension_setting} 
\parbox{1.1in}{
\begin{tikzpicture}[description/.style={fill=white,inner sep=2pt}]
\matrix (m) [ampersand replacement=\&, matrix of math nodes, row sep=2em, column sep=2em, text height=1.5ex, text depth=0.5ex]
{ 
  \C \& {} \& \E \\
  {} \& \D \\ 
};
\path[->] 
(m-1-1) edge node[above] {\footnotesize$F$} (m-1-3) edge node[below,sloped] {\footnotesize$K$} (m-2-2)
;
\end{tikzpicture}
}
\end{equation} 
and also that we want to extend $F$ along $K$, that is, to construct a functor $G \colon \D \too \E$ such that $F = G \circ K$. Some problems may arise regarding this matter: 
\begin{itemize}[itemsep=2pt,topsep=0pt]
\item[$\bullet$] Two morphisms in $\C$ may have different images under $F$, but equal under $K$.

\item[$\bullet$] There may exist $X, Y \in \C$ such that $\Hom_{\C}(X,Y) = \emptyset$, $\Hom_{\E}(F(X),F(Y)) = \emptyset$,  and  $\Hom_{\D}(K(X),K(Y)) \neq \emptyset$.
\end{itemize}
Due to these inconveniences, it seems more reasonable to extend $F$ along $K$ by constructing a functor $G \colon \D \too \E$ along with a natural transformation $F \Rightarrow G \circ K$. Here is where the concept of Kan extensions appears.

\begin{definition}\label{def:ext_Kan}
Given two functors $F \colon \C \too \E$ and $K \colon \C \too \D$, a (\textbf{left}) \textbf{Kan extension of $\bm{F}$ along $\bm{K}$} is a functor $G \colon \D \too \E$, along with a natural transformation $\eta \colon F \Rightarrow G \circ K$, satisfying the following:
\begin{itemize}
[itemsep=2pt,topsep=0pt]
\item[$\bullet$] \textbf{\underline{Universal property}:} For any other pair $(G' \colon \D \too \E, \eta' \colon F \Rightarrow G' \circ K)$, there exists a unique natural transformation $\alpha \colon G \Rightarrow G'$ such that $\eta' = \alpha_K \circ \eta$. 
\end{itemize}
\end{definition}

Graphically, we have the following diagram of functors and natural transformations:
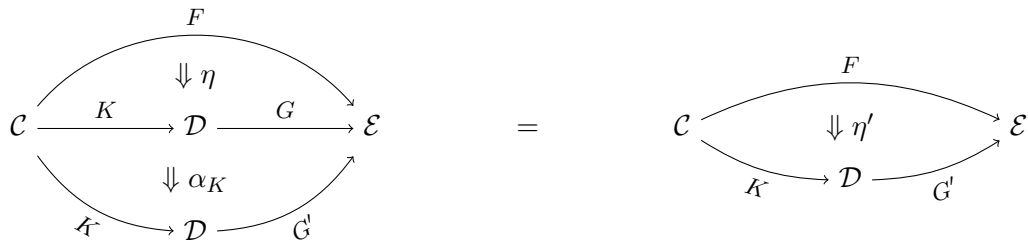
\begin{figure}[H]
\centering
\begin{tikzpicture}[description/.style={fill=white,inner sep=2pt}]
\matrix (m) [matrix of math nodes, row sep=-0.5em, column sep=4em, text height=2.25ex, text depth=1.25ex]
{ 
  {} & \Downarrow\eta & {} & {} & {} & {} \\
  \C & \D & \E & = & \C & \Downarrow\eta' & \E \\
  {} & \Downarrow \alpha_K & {} & {} & {} & \D \\
  {} & \D & {} & {} & {} & {} \\
};
\path[->]
(m-2-1) edge node[above] {\footnotesize$K$} (m-2-2)
(m-2-2) edge node[above] {\footnotesize$G$} (m-2-3)
(m-2-1) edge [bend left=50] node[above] {\footnotesize$F$} (m-2-3)
(m-2-1) edge [bend right=25] node[below,sloped] {\footnotesize$K$} (m-4-2)
(m-4-2) edge [bend right=25] node[below,sloped] {\footnotesize$G'$} (m-2-3)
(m-2-5) edge [bend right=15] node[below,sloped] {\footnotesize$K$} (m-3-6)
(m-3-6) edge [bend right=15] node[below,sloped] {\footnotesize$G'$} (m-2-7)
(m-2-5) edge [bend left=25] node[above] {\footnotesize$F$} (m-2-7)
;
\end{tikzpicture}
\caption[Universal property of Kan extensions]{Universal property of Kan extensions.}
\end{figure}

\newpage

Note that Kan extensions are unique up to natural isomorphisms. That is why we will refer to $G$ as \emph{the} Kan extension of $F$ along $K$, and will be denoted by:
\[
G = \Lan_K(F).
\] 
The notation $\Lan$ is simply a contraction of the term ``left Kan''.

Presenting most of the interesting examples of Kan extensions would require to recall several concepts from algebraic topology, model category theory, and graph theory. For this reason, examples will be omitted. However, we suggest the reader to check the book \emph{\cite[Chapter 1]{RiehlBook}} by E. Riehl. For examples more focused in algebraic topology, the notes  \emph{\cite{Riehl}}, also by Riehl, are an excellent source. For examples in model category theory, we suggest P. S. Hirschhorn's \emph{\cite[Proposition 8.4.4]{Hir}}. Finally, for the reader interested in vector spaces or directed graph, we recommend to check M. C. Lehner (B. Sc.) thesis \emph{\cite{Lehner}}. One can also find examples related to database migration, developed by D. I. Spivak in \emph{\cite{SpivakMigration}}, which are also presented in Spivak's book \emph{\cite[Section 7.1.4.6]{SpivakBook}} in a more expository way.


\subsection*{Kan extensions and representable functors}

Consider a functor $K \colon \C \too \D$. For each category $\E$, we have that $K$ induces a contravariant functor 
\[
- \circ K \colon \Fun(\D,\E) \too \Fun(\C,\E)
\]
given by $(- \circ K)(G) := G \circ K$, for every $G \in \Fun(\D,\E)$. To define $- \circ K$ on morphisms in $\Fun(\D,\E)$, suppose we are given a natural transformation $\omega \colon G \Rightarrow G'$, with $G, G' \in \Fun(\D,\E)$, then define $\omega \circ K \colon G \circ K \Rightarrow G' \circ K$ that the natural transformation $\omega_K$. 

This functor $- \circ K$ helps us to note that the Kan extension of a functor $F \colon \C \too \E$ along $K$ is a \emph{representation} for the functor 
\[
\Hom_{\Fun(\C,\E)}(F, - \circ K) \colon \Fun(\D,\E) \too \Set,
\]
that is, there exists a natural isomorphism 
\[
\Hom_{\Fun(\C,\E)}(F, - \circ K) \cong \Hom_{\Fun(\D,\E)}(\Lan_K(F),-)
\]
given by $\eta$. Indeed, for each $G \in \Fun(\D,\E)$, we have a bijection of sets
\begin{align*}
\eta_G \colon \Hom_{\Fun(\D,\E)}(\Lan_K(F),G) & \xrightarrow{\sim} \Hom_{\Fun(\C,\E)}(F,G \circ K) \\
\alpha & \mapsto \alpha_K \circ \eta.
\end{align*}
Moreover, for each morphism $\omega \colon G \Rightarrow G'$ in $\Fun(\D,\E)$, we have the following commutative diagram:
\[
\begin{tikzpicture}[description/.style={fill=white,inner sep=2pt}]
\matrix (m) [matrix of math nodes, row sep=3em, column sep=3em, text height=2.25ex, text depth=1.25ex]
{ 
  \Hom_{\Fun(\D,\E)}(\Lan_K(F), G) & \Hom_{\Fun(\C,\E)}(F, G \circ K) \\
  \Hom_{\Fun(\D,\E)}(\Lan_K(F), G') & \Hom_{\Fun(\C,\E)}(F, G' \circ K) \\
};
\path[->]
(m-1-1) edge node[above] {\footnotesize$\sim$} node[below] {\footnotesize$\eta_G$} (m-1-2) edge node[left] {\footnotesize$\Hom_{\Fun(\D,\E)}(\Lan_K(F),\omega)$} (m-2-1)
(m-1-2) edge node[right] {\footnotesize$\Hom_{\Fun(\C,\E)}(F, \omega_K)$} (m-2-2)
(m-2-1) edge node[above] {\footnotesize$\sim$} node[below] {\footnotesize$\eta_{G'}$} (m-2-2)
;
\end{tikzpicture}
\]
By ``commutative'' we mean that the equality of functions 
\[
\Hom_{\mathsf{Fun}(\C,\E)}(F,\omega_K) \circ \eta_G = \eta_{G'} \circ \Hom_{\mathsf{Fun}(\D,\E)}(\text{Lan}_K(F),\omega)
\]
holds, which is easy to verify.


\subsection*{Kan extensions as adjunctions}

In this section we present a first contact to the interplay between Kan extensions and adjoint functors. Suppose we are given functors $F$ and $K$ as in \eqref{fig:extension_setting}. Under certain conditions, the functor $- \circ K \colon \Fun(\D,\E) \too \Fun(\C,\E)$ introduced previously, has a left adjoint. Namely, assuming that the Kan extension of $F$ along $K$ exists for every $F \in \Fun(\C,\E)$ (for example, this occurs in the case $\C$ is small, or skeletally small, and $\E$ is cocomplete, studied in \emph{Section~\ref{sec:formula}}), we have a functor 
\[
\Lan_K(-) \colon \Fun(\C,\E) \too \Fun(\D,\E)
\] 
defined as $F \mapsto \Lan_K(F)$ on objects of $\Fun(\C,\E)$. In order to define $\Lan_K(-)$ on morphisms of $\Fun(\C,\E)$, we need to use the universal property of Kan extensions. Suppose we are given a morphism $\nu \colon F \Rightarrow F'$ in $\Fun(\C,\E)$. For $\Lan_K(F)$ and $\Lan_K(F')$ we have natural transformations $\eta \colon F \Rightarrow \Lan_K(F) \circ K$ and $\eta' \colon F' \rightarrow \Lan_K(F') \circ K$ which satisfy the universal property of Kan extensions, and so there exists a unique natural transformation 
\[
\Lan_K(\nu) \colon \Lan_K(F) \Rightarrow \Lan_K(F')
\] 
such that $\eta' \circ \nu = \Lan_K(\nu)_K \circ \eta$. It follows that the mapping $\nu \mapsto \text{Lan}_K(\nu)$ is well defined. Using the universal property again, one can also show that 
\[
\Lan_K(\nu' \circ \nu) = \Lan_K(\nu') \circ \Lan_K(\nu) \mbox{ \ and \ } \Lan_K({\rm id}_F) = {\rm id}_{\Lan_K(F)},
\] 
for each pair of natural transformations $\nu \colon F \Rightarrow F'$ and $\nu \colon F' \Rightarrow F''$. 

On the other hand, we have a natural isomorphism 
\begin{align*}
\eta_G \colon \Hom_{\Fun(\D,\E)}(\Lan_K(F), G) & \cong \Hom_{\Fun(\C,\E)}(F, G \circ K),
\end{align*}
implying that $(\Lan_K(-), - \circ K)$ is an \textbf{adjoint pair}, that is, $\Lan_K(-)$ is a left adjoint of $- \circ K$, denoted $\Lan_K(-) \dashv (- \circ K)$.
\[
\begin{tikzpicture}[description/.style={fill=white,inner sep=2pt}]
\matrix (m) [matrix of math nodes, row sep=0em, column sep=1.5em, text height=2.25ex, text depth=1.25ex]
{ 
  \Fun(\C,\E) & \perp & \Fun(\D,\E) \\
};
\path[->]
(m-1-3) edge [bend left=15] node[below] {\footnotesize$- \circ K$} (m-1-1)
(m-1-1) edge [bend left=15] node[above] {\footnotesize$\Lan_K(-)$} (m-1-3)
;
\end{tikzpicture}
\]


\subsection*{Colimits as Kan extensions}\label{sec:universal}

The following quote is due to Saunders Mac Lane, and can be found in his book \emph{\cite[Chapter 10]{MacLane}}:
\[
\mbox{``\emph{The notion of Kan extensions subsumes all the other fundamental concepts of Category Theory}''.}
\]
Basically, what this quote suggest is that Kan extensions are the most universal constructions in category theory. In this section, we want to emphasize this by showing that the concept of colimit is a particular Kan extension. Moreover, we complement Mac Lane's quote by proving later that adjoint pairs and Kan extensions are somehow equivalent. 

Denote by $\tau$ the \textbf{terminal category}, that is, $\tau$ has only one object, say $\ast$, and a unique morphism $\ast \to \ast$, namely, the identity on $\ast$. Note that $\tau$ is the terminal object in $\Cat$, the category of small categories and functors, thus suggesting the name ``terminal category'' for $\tau$.

\begin{proposition}[colimits vs. Kan extensions]
The colimit of a functor $F \colon \C \too \E$ exists if, and only if, its Kan extension along the only functor $K \colon \C \too \tau$ also exists.
\end{proposition}

\begin{myproof}
We only prove the ``only if'' part, since the ``if'' part follows in a similar way.

First note that $K \colon \C \too \tau$ is the constant functor given by $C \mapsto \ast$, for every $C \in \C$, and by $f \mapsto {\rm id}_{\ast}$ for every $f \in \Mor(\C)$. Let $X \in \E$ be the colimit of $F$. Then, we have a family of commutative triangles 
\[
\begin{tikzpicture}[description/.style={fill=white,inner sep=2pt}]
\matrix (m) [matrix of math nodes, row sep=1.5em, column sep=1.5em, text height=2.25ex, text depth=1.25ex]
{ 
  F(C) & {} & F(C') \\
  {} & X & {} \\
};
\path[->]
(m-1-1) edge node[above] {\footnotesize$F(f)$} (m-1-3) edge node[below,sloped] {\footnotesize$\eta_C$} (m-2-2)
(m-1-3) edge node[below,sloped] {\footnotesize$\eta_{C'}$} (m-2-2)
;
\end{tikzpicture}
\]  
with $f$ running over $\Hom_{\C}(C,C')$, satisfying the universal property of colimits. Define the functor $G \colon \tau \too \E$ by $G(\ast) := X$ and $G({\rm id}_{\ast}) := {\rm id}_X$. Let us verify that $G = \Lan_K(F)$. First, note that the family  $\{ \eta_C \colon F(C) \to X \}_{C \in \C}$ defines a natural transformation $\eta \colon F \Rightarrow G  \circ K$. 

We now check that the pair $(G,\eta)$ satisfies the universal property of Kan extensions. Suppose we have another functor $G' \colon \tau \too \E$ along with a natural transformation $\eta' \colon F \Rightarrow G' \circ K$. We construct a natural transformation $\alpha \colon G \Rightarrow G'$ such that $\eta' = \alpha_K \circ \eta$. Note that $\alpha$ is going to be determined by a unique morphism $\alpha_\ast \colon X \to G'(\ast)$. Since $\eta' \colon F \Rightarrow G' \circ K$ is a natural transformation, one has that $\eta'_{C'} \circ F(f) = \eta'_C$ for every $f \in \Hom_{\C}(C,C')$. By the universal property of colimits, we can find such $\alpha_\ast$ satisfying $\alpha_\ast \circ \eta_C = \eta'_C$, for every $C \in \C$. This family of equalities can be represented as $\eta' = \alpha_K \circ \eta$. Finally, the fact that $\alpha$ is the only natural transformation $G \Rightarrow G'$ satisfying the previous equality, follows by a straightforward application of the universal property of colimits. 
\end{myproof}


\subsection*{Adjoint pairs vs. Kan extensions}

In this section we study the interaction between Kan extensions and adjoint pairs, one of the most important notions in category theory. Let us first see how Kan extensions can be obtained from an adjoint pair.

\begin{proposition}[from adjoint pairs to Kan extensions]\label{prop:de_pares_adjuntos_a_extensiones_de_Kan}
Let $F \colon \C \too \D$ and $G \colon \D \too \C$ be two functors such that $(F,G)$ is an adjoint pair. Then, $G = \Lan_F({\rm id}_{\C})$. In this case, the natural transformation $\eta \colon {\rm id}_{\C} \Rightarrow G \circ F$ satisfying \emph{Definition~\ref{def:ext_Kan}} coincides with the unit of $(F,G)$.
\end{proposition}

\begin{myproof}
Let $\eta \colon {\rm id}_{\C} \Rightarrow G \circ F$ be the unit of the adjunction $F \dashv G$. Let us prove that $(G,\eta)$ satisfies the universal property of Kan extensions for ${\rm id}_{\C}$ along $F$. Suppose we are given another functor $G' \colon \D \too \C$ along with a natural transformation $\eta' \colon {\rm id}_{\C} \Rightarrow G' \circ F$. Consider $G(D) \in \C$. Using the counit $\epsilon \colon F \circ G \Rightarrow {\rm id}_{\D}$ of $(F,G)$ and the natural transformation $\eta'$, we have two morphisms $\eta'_{G(D)} \colon G(D) \to G' \circ F \circ G(D)$ and $\epsilon_D \colon F \circ G(D) \to D$, from which we define $\alpha_D$ as:
\[
\alpha_D := G'(\epsilon_D) \circ \eta'_{G(D)}, \mbox{ for every $D \in \D$}.
\]
Now for every $g \in \Hom_{\D}(D,D')$, using the fact that $\epsilon$ and $\eta'$ are natural transformations, one can verify that $G'(g) \circ \alpha_D = \alpha_{D'} \circ G(g)$, that is, the family of morphisms $\alpha := \{ \alpha_D \colon G(D) \to G'(D) \}_{D \in \D}$ defines a natural transformation $\alpha \colon G \Rightarrow G'$. We show that such $\alpha$ satisfies the conditions in \emph{Definition~\ref{def:ext_Kan}}. 

Before showing the equality $\eta' = \alpha_F \circ \eta$, recall that the unit $\eta \colon {\rm id}_{\C} \Rightarrow G \circ F$ and counit $\epsilon \colon F \circ G \Rightarrow {\rm id}_{\D}$ of the adjoint pair $(F,G)$ satisfy the equalities
\begin{align}
\epsilon_{F(C)} \circ F(\eta_C) & = {\rm id}_{F(C)}, \mbox{ for every } C \in \C, \label{eqn:triangular_uno} \\
G(\epsilon_D) \circ \eta_{G(D)} & = {\rm id}_{G(D)}, \mbox{ for every } D \in \D, \label{eqn:triangular_dos}
\end{align}
which are known as the \textbf{triangle identities} (See \emph{\cite[Lemma 2.2.2]{Leinster}}). Now, using the fact that $\eta'$ is a natural transformation, along with $\eqref{eqn:triangular_uno}$, we have that $(\alpha_F \circ \eta)_C = \eta'_C$, for every $C \in \C$. Hence, $\eta' = \alpha_F \circ \eta$.

Finally, suppose that there exists another natural transformation $\alpha' \colon G \Rightarrow G'$ such that $\eta' = \alpha'_F \circ \eta$. Let us show that $\alpha = \alpha'$. For every $D \in \D$, we have:
\begin{align*}
\alpha'_D & = \alpha'_D \circ G(\epsilon_D) \circ \eta_{G(D)} & \mbox{(by \eqref{eqn:triangular_dos})} \\
& = G'(\epsilon_D) \circ \alpha'_{F(G(D))} \circ \eta_{G(D)} & \mbox{(since $\alpha'$ is a natural transformation)} \\
& = G'(\epsilon_D) \circ \alpha_{F(G(D))} \circ \eta_{G(D)} & \mbox{(since $\eta' = \alpha'_F \circ \eta$)} \\
& = \alpha_D \circ G(\epsilon_D) \circ \eta_{G(D)} & \mbox{(since $\alpha$ is a natural transformation)} \\
& = \alpha_D. & \mbox{(by \emph{\eqref{eqn:triangular_dos}})}
\end{align*}
Therefore, $G = \Lan_F({\rm id}_{\C})$, along with the natural transformation $\eta \colon {\rm id}_{\C} \Rightarrow G \circ F$. 
\end{myproof}

A natural question arising at this point is under which condition it is possible to show the converse of \emph{Proposition~\ref{prop:de_pares_adjuntos_a_extensiones_de_Kan}}. That is, how can we obtain adjoint pairs from Kan extensions of ${\rm id}_{\C}$? The condition we are interested in is associated to the concept of functors that preserve Kan extensions, stated below.

\newpage

\begin{definition}
Let $F \colon \C \too \E$ be a functor with Kan extension along $K \colon \C \too \D$, say $(\Lan_K(F), \eta)$. A functor $H \colon \E \too \F$ \textbf{preserves} $(\Lan_K(F), \eta)$ if the pair $(H \circ \Lan_K(F), H(\eta))$ is the Kan extension of $H \circ F$ along $K$.
\begin{figure}[H]
\centering
\begin{tikzpicture}[description/.style={fill=white,inner sep=2pt}]
\matrix (m) [matrix of math nodes, row sep=3em, column sep=4em, text height=2.25ex, text depth=1.25ex]
{ 
  \C & \E & \F \\
  {} & \D & {} \\
};
\path[->]
(m-1-2)-- node[pos=0.3] {\footnotesize$H(\eta)\Downarrow$} (m-2-2)
(m-1-1) edge node[above] {\footnotesize$F$} (m-1-2) edge node[below,sloped] {\footnotesize$K$} (m-2-2)
(m-1-2) edge node[above] {\footnotesize$H$} (m-1-3)
(m-2-2) edge node[right] {\footnotesize$\mbox{ \ }H \circ \Lan_K(F) = {\rm Lan}_K(H \circ F)$} (m-1-3)
;
\end{tikzpicture}
\caption[Preservation of Kan extensions]{Preservation of Kan extensions.}
\end{figure}
\end{definition}

The following result is the dual of \emph{\cite[Theorem 4.5]{Lehner}}.

\begin{proposition}[from Kan extensions to adjoint pairs]\label{prop:from_Kan_extensions_to_adjoint_pairs}
If $(G,\eta)$ is the Kan extension of ${\rm id}_{\C}$ along $F$, and if $F$ preserves $(G,\eta)$, then $(F,G)$ is an adjoint pair whose unit is given by $\eta$. 
\end{proposition}

Before giving the proof of the previous proposition, we show the following property of adjoint pairs, whose proof can be found also in  \emph{\cite[Theorem 4.2]{Lehner}}, where the author uses different arguments.

\begin{proposition}
Left adjoint functors preserve (left) Kan extensions.
\end{proposition}

\begin{myproof}
Suppose we are given the following diagram of functors:
\[
\begin{tikzpicture}[description/.style={fill=white,inner sep=2pt}]
\matrix (m) [matrix of math nodes, row sep=2em, column sep=2em, text height=1.25ex, text depth=0.25ex]
{ 
  \C & {} & \E & \perp & \F \\
  {} & \D & {} & {} & {} \\
};
\path[->]
(m-1-1)-- node[pos=0.4] {\footnotesize$\mbox{ \ \ \ \ }\eta\Downarrow$} (m-2-3)
(m-1-1) edge node[above] {\footnotesize$F$} (m-1-3) edge node[below,sloped] {\footnotesize$K$} (m-2-2)
(m-1-3) edge [bend left=20] node[above] {\footnotesize$H$} (m-1-5)
(m-1-5) edge [bend left=20] node[below] {\footnotesize$Q$} (m-1-3)
(m-2-2) edge node[below,sloped] {\footnotesize$G$} (m-1-3)
;
\end{tikzpicture}
\]
where $G = \Lan_K(F)$ and $H$ is a left adjoint of $Q$. 

Let us show $H \circ G = \Lan_K(H \circ F)$. Assume that we are given a functor $G' \colon \D \too \F$ along with a natural transformation $\gamma \colon H \circ F \Rightarrow G' \circ K$. 

From the natural transformations $Q(\gamma) \colon Q \circ H \circ F \Rightarrow Q \circ G' \circ K$ and $\mu_F \colon F \Rightarrow Q \circ H \circ F$, where $\mu \colon {\rm id}_{\E} \Rightarrow Q \circ H$ is the unit of the adjoint pair $(H, Q)$, we define
\[
\eta' := Q(\gamma) \circ \mu_F \colon F \Rightarrow Q \circ G' \circ K.
\]
Since $G = \Lan_K(F)$, there exists a unique natural transformation $\alpha^0 \colon G \Rightarrow Q \circ G'$ such that $\eta' = \alpha^0_K \circ \eta$. Now, considering also the counit $\nu \colon H \circ Q \Rightarrow {\rm id}_{\mathcal{F}}$ of $(H,Q)$, we have natural transformations $\nu_{G'} \colon H \circ Q \circ G' \Rightarrow G'$ and $H(\alpha^0) \colon H \circ G \Rightarrow H \circ Q \circ G'$, from which we define:
\[
\alpha := \nu_{G'} \circ H(\alpha^0) \colon H \circ G \Rightarrow G'.
\]
We show that such $\alpha$ satisfies the conditions in the universal property of \emph{Definition~\ref{def:ext_Kan}}. 

For every $C \in \C$, we have:
\begin{align*}
(\alpha_K \circ H(\eta))_C & = \alpha_{K(C)} \circ H(\eta_C) \\
& = \nu_{G' (K(C))} \circ H(\alpha^0_{K(C)} \circ \eta_C) \\
& = \nu_{G' (K (C))} \circ H(\eta'_C) \\
& = \nu_{G' (K(C))} \circ H(Q(\gamma_C)) \circ H(\mu_{F(C)}) \\
& = \gamma_C \circ \nu_{H (F(C))} \circ H(\mu_{F(C)}) & \mbox{(since $\nu$ is a natural transformation)} \\
& = \gamma_C. & \mbox{(by \emph{\eqref{eqn:triangular_uno}})} 
\end{align*}
Then, we have $\alpha_K \circ H(\eta) = \gamma$. Finally, suppose there is another natural transformation $\alpha' \colon H \circ G \Rightarrow G'$ such that $\alpha'_K \circ H(\eta) = \gamma$. Then, from $Q(\alpha') \colon Q \circ H \circ G \Rightarrow Q \circ G'$ and $\mu_G \colon G \Rightarrow Q \circ H \circ G$, we define
\[
\overline{\alpha}^0 := Q(\alpha') \circ \mu_G \colon G \Rightarrow Q \circ G'.
\]
Using the fact that $\mu$ is a natural transformation, we have $\overline{\alpha}^0_K \circ \eta = \eta'$. Hence, $\overline{\alpha}^0 = \alpha^0$, since $G = \Lan_K(F)$. From this fact, we can see that $\alpha = \alpha'$. Indeed, for every $D \in \D$, we have:
\begin{align*}
\alpha_D & = \nu_{G'(D)} \circ H(\alpha^0_D) = \nu_{G'(D)} \circ H(\overline{\alpha}^0_D) \\
& = \nu_{G'(D)} \circ H(Q(\alpha'_D) \circ \mu_{G(D)}) \\
& = \nu_{G'(D)} \circ H(Q(\alpha'_D)) \circ H(\mu_{G(D)}) \\
& = \alpha'_D \circ \nu_{H(G(D))} \circ H(\mu_{G(D)}) & \mbox{(since $\nu$ is a natural transformation)} \\
& = \alpha'_D \circ {\rm id}_{H(G(D))} = \alpha'_D. & \mbox{(by \emph{\eqref{eqn:triangular_uno}})} 
\end{align*}
Therefore, $\alpha = \alpha'$, and thus $H \circ G = \Lan_K(H \circ F)$. 
\end{myproof}

\begin{corollary}\label{coro:adjunciones_y_extensiones_de_Kan_inducidas}
If $(F,G)$ is an adjoint pair with unit $\eta$, then $F$ preserves $(G,\eta)$, where $(G,\eta)$ is the Kan extension given by \emph{Proposition~\ref{prop:de_pares_adjuntos_a_extensiones_de_Kan}}.
\end{corollary}

We are now ready to give a proof of \emph{Proposition~\ref{prop:from_Kan_extensions_to_adjoint_pairs}}. However, we present a more complete result.

\begin{theorem}[adjoint pairs vs. Kan extensions]\label{theo:Kan_extensions_vs_adjoint_pairs}
A functor $F \colon \C \too \D$ has a right adjoint $G \colon \D \too \C$ if, and only if, $G = \Lan_F({\rm id}_{\C})$ and $\Lan_F({\rm id}_{\C})$ is preserved by $F$. In this case, $F \dashv \Lan_F({\rm id}_{\C})$ and $\eta \colon {\rm id}_{\C} \Rightarrow \Lan_F({\rm id}_{\C}) \circ F$ is the unit of the adjunction.  
\end{theorem}

\begin{myproof}
The ``only if'' part follows from \emph{Proposition~\ref{prop:de_pares_adjuntos_a_extensiones_de_Kan}} and \emph{Corollary~\ref{coro:adjunciones_y_extensiones_de_Kan_inducidas}}.

Now suppose that ${\rm id}_{\C}$ has Kan extension along $F$, say $(G,\eta)$, and that it is preserved by $F$. Let us see that $F \dashv G$. By  \emph{\cite[Corollary 2.2.6]{Leinster}}, it suffices to obtain natural transformations $\eta \colon {\rm id}_{\C} \Rightarrow G \circ F$ and $\epsilon \colon F \circ G \Rightarrow {\rm id}_{\D}$ satisfying the triangle identities:   
\begin{itemize}[itemsep=2pt,topsep=0pt]
\item[$\iroman$] ${\rm id}_F = \epsilon_F \circ F(\eta)$.

\item[$\iiroman$] ${\rm id}_G = G(\epsilon) \circ \eta_G$.
\end{itemize}
Let $\eta \colon {\rm id}_{\C} \Rightarrow G \circ F$ be the natural transformation accompanying the Kan extension $G$. We know that $F$ preserves $(G,\eta)$, that is, $F \circ G = \Lan_F(F)$. Then, there exists a unique natural transforamtion $\epsilon \colon F \circ G \Rightarrow {\rm id}_{\D}$ such that ${\rm id}_F = \epsilon_F \circ F(\eta)$, which proves $\iroman$.
\[
\begin{tikzpicture}[description/.style={fill=white,inner sep=2pt}]
\matrix (m) [matrix of math nodes, row sep=-0.5em, column sep=4em, text height=2.25ex, text depth=1.25ex]
{ 
  {} & \Downarrow F(\eta) & {} & {} & {} & {} \\
  \C & \D & \D & = & \C & \Downarrow {\rm id}_F & \D \\
  {} & \Downarrow \epsilon_F & {} & {} & {} & \D \\
  {} & \D & {} & {} & {} & {} \\
};
\path[->]
(m-2-1) edge node[above] {\footnotesize$F$} (m-2-2)
(m-2-2) edge node[above] {\footnotesize$F \circ G$} (m-2-3)
(m-2-1) edge [bend left=50] node[above] {\footnotesize$F$} (m-2-3)
(m-2-1) edge [bend right=25] node[below,sloped] {\footnotesize$F$} (m-4-2)
(m-4-2) edge [bend right=25] node[below,sloped] {\footnotesize${\rm id}_{\D}$} (m-2-3)
(m-2-5) edge [bend right=15] node[below,sloped] {\footnotesize$F$} (m-3-6)
(m-3-6) edge [bend right=15] node[below,sloped] {\footnotesize${\rm id}_{\D}$} (m-2-7)
(m-2-5) edge [bend left=25] node[above] {\footnotesize$F$} (m-2-7)
;
\end{tikzpicture}
\]
We can obtain $\iiroman$ by showing the equality $(G(\epsilon) \circ \eta_G)_F \circ \eta = \eta$, using the fact that $\eta$ is a natural transformation, along with $\iroman$. This implies that $G(\epsilon) \circ \eta_G = {\rm id}_{G}$, by the universal property of $(G,\eta)$.
\[
\begin{tikzpicture}[description/.style={fill=white,inner sep=2pt}]
\matrix (m) [matrix of math nodes, row sep=-0.5em, column sep=3em, text height=2.25ex, text depth=1.25ex]
{ 
  {} & \Downarrow \eta & {} & {} & {} & {} \\
  \C & \D & \C & = & \C & \Downarrow \eta & \C \\
  {} & \Downarrow {\rm id}_{G \circ F} \Downarrow (G(\epsilon) \circ \eta_G)_F & {} & {} & {} & \D \\
  {} & \D & {} & {} & {} & {} \\
};
\path[->]
(m-2-1) edge node[above] {\footnotesize$F$} (m-2-2)
(m-2-2) edge node[above] {\footnotesize$G$} (m-2-3)
(m-2-1) edge [bend left=35] node[above] {\footnotesize${\rm id}_{\C}$} (m-2-3)
(m-2-1) edge [bend right=20] node[below,sloped] {\footnotesize$F$} (m-4-2)
(m-4-2) edge [bend right=20] node[below,sloped] {\footnotesize$G$} (m-2-3)
(m-2-5) edge [bend right=15] node[below,sloped] {\footnotesize$F$} (m-3-6)
(m-3-6) edge [bend right=15] node[below,sloped] {\footnotesize$G$} (m-2-7)
(m-2-5) edge [bend left=30] node[above] {\footnotesize${\rm id}_{\C}$} (m-2-7)
;
\end{tikzpicture}
\] 
\end{myproof}


\section{Coends to compute Kan extensions}\label{sec:formula}

Under certain conditions on $\C$, $\D$ and $\E$ in \eqref{fig:extension_setting}, it is possible to compute Kan extensions by the use of a formula which involves certain colimit, called \emph{coend}. Below we state such a formula, for which we give a detailed proof. In \emph{\cite[Chapter 10, Section 4, Theorem 1]{MacLane}}, the reader can find a different proof written by S. Mac Lane.

\newpage

\begin{theorem}[Existence of Kan extensions from coends]\label{theo:MacLane}
Suppose that $\C$ is a small category and that $\E$ is a cocomplete category, that is, $\E$ has colimits. Then, given two functors $F \colon \C \too \E$ and $K \colon \C \too \D$, the Kan extension of $F$ along $K$ exists, and it is computed by the expression  
\begin{align}\label{eqn:Grothendieck}
\Lan_K(F)(D) & = \displaystyle\operatorname*{\int}^{C \in \C} \Hom_{\D}(K(C),D) \odot F(C), 
\end{align}
for every $D \in \D$.
\end{theorem}

The proof of the previous theorem is rather technical, although not difficult. In order to understand the meaning of the previous formula, it is necessary to recall the displayed notation.


\subsection*{The copower bifunctor}

The symbol $\odot$ in \emph{\eqref{eqn:Grothendieck}} is called \textbf{copower}, and it is defined as a bifunctor $\odot \colon \Set \times \E \too \E$ as follows: if $S \in \Set$ and $E \in \E$, then $S \odot E$ is the coproduct of copies of $E$ indexed over the set $S$. The natural inclusions of this coproduct will be denoted by $i^E_s \colon E \hookrightarrow S \odot E$, for each $s \in S$. 

Now let us see how to define $- \odot -$ on morphisms of $\Set$ and $\E$. First, suppose we are given a function $t \colon S \to U$. By the universal property of coproducts, there exists a unique morphism $t \odot E$ in $\E$ such that the following diagram commutes:
\begin{equation}\label{fig:t_tensor_E} 
\parbox{1.1in}{
\begin{tikzpicture}[description/.style={fill=white,inner sep=2pt}]
\matrix (m) [ampersand replacement=\&, matrix of math nodes, row sep=2em, column sep=4em, text height=1.75ex, text depth=0.25ex]
{ 
  E \& S \odot E \\
  {} \& U \odot E \\
};
\path[->]
(m-1-1) edge node[above] {\footnotesize$i^E_s$} (m-1-2)
(m-1-1) edge node[below,sloped] {\footnotesize$i^E_{t(s)}$} (m-2-2)
;
\path[dotted,->]
(m-1-2) edge node[left] {\footnotesize$\exists!$} node[right] {\footnotesize$t \odot E$} (m-2-2)
;
\end{tikzpicture}
}
\end{equation}
On the other hand, suppose we have a morphism $f \colon E_1 \to E_2$ in $\E$. Using again the universal property of coproducts, we have that there exists a unique morphism $S \odot f$ in $\E$ such that the following diagram commutes: 
\begin{equation}\label{fig:S_tensor_f} 
\parbox{1.1in}{
\begin{tikzpicture}[description/.style={fill=white,inner sep=2pt}]
\matrix (m) [ampersand replacement=\&, matrix of math nodes, row sep=2em, column sep=4em, text height=1.75ex, text depth=0.25ex]
{ 
  E_1 \& S \odot E_1 \\
  E_2 \& S \odot E_2 \\
};
\path[->]
(m-1-1) edge node[above] {\footnotesize$i^{E_1}_{s}$} (m-1-2)
(m-2-1) edge node[below] {\footnotesize$i^{E_2}_{s}$} (m-2-2)
(m-1-1) edge node[left] {\footnotesize$f$} (m-2-1)
;
\path[dotted,->]
(m-1-2) edge node[left] {\footnotesize$\exists!$} node[right] {\footnotesize$S \odot f$} (m-2-2)
;
\end{tikzpicture}
}
\end{equation}
One can also verify that $\odot \colon \Set \times \E \too \E$, given by the mappings $(S,E) \mapsto E^{(S)}$, $(t,E) \mapsto t \odot E$, and $(S,f) \mapsto S \odot f$, is indeed a bifunctor. The details are left to the reader.


\subsection*{The coend of a bifunctor}

The symbol $\int$ is called \emph{coend}, and is defined below.

\begin{definition}\label{def:coend}\label{def:cofin}
Let $H \colon \C\op \times \C \too \E$ be a bifunctor. A \textbf{coend} of $H$ is an object $E \in \E$, along with a family of morphisms $\lambda_C \colon H(C,C) \to E$ in $\E$, with $C$ running over $\C$, such that the following two properties hold:
\begin{itemize}[itemsep=2pt,topsep=0pt]
\item[$\bullet$] \textbf{\underline{Naturality}:} For every $f \in \Hom_{\C}(C,C')$, the following diagram in $\E$ commutes:
\[
\begin{tikzpicture}[description/.style={fill=white,inner sep=2pt}]
\matrix (m) [matrix of math nodes, row sep=3em, column sep=5em, text height=2ex, text depth=0.25ex]
{ 
  H(C',C) & H(C',C') \\
  H(C,C) & E \\
};
\path[->]
(m-1-1) edge node[above] {\footnotesize$H(C',f)$} (m-1-2)
(m-2-1) edge node[below] {\footnotesize$\lambda_C$} (m-2-2)
(m-1-1) edge node[left] {\footnotesize$H(f,C)$} (m-2-1)
(m-1-2) edge node[right] {\footnotesize$\lambda_{C'}$} (m-2-2)
;
\end{tikzpicture}
\]
That is, $\lambda_C \circ H(f,C) = \lambda_{C'} \circ H(C',f)$.

\item[$\bullet$] \textbf{\underline{Universal property}:} Given another object $E' \in \E$, along with a family of natural morphisms $\lambda'_C \colon H(C,C) \to E'$ in $\E$, with $C$ running over $\C$, that is, such that $\lambda'_C \circ H(f,C) = \lambda'_{C'} \circ H(C',f)$ for every $f \in \Hom_{\C}(C,C')$, there exists a unique morphism $\lambda \colon E \to E'$ in $\E$ such that $\lambda \circ \lambda_C = \lambda'_C$, for every $C \in \C$. 
\begin{figure}[H]
\centering
\begin{tikzpicture}[description/.style={fill=white,inner sep=2pt}]
\matrix (m) [matrix of math nodes, row sep=3em, column sep=5em, text height=2ex, text depth=0.25ex]
{ 
  H(C',C) & H(C',C') \\
  H(C,C) & E \\
  {} & {} & E' \\
};
\path[->]
(m-1-1) edge node[above] {\footnotesize$H(C',f)$} (m-1-2)
(m-2-1) edge node[below] {\footnotesize$\lambda_C$} (m-2-2)
(m-1-1) edge node[left] {\footnotesize$H(f,C)$} (m-2-1)
(m-1-2) edge node[right] {\footnotesize$\lambda_{C'}$} (m-2-2)
(m-1-2) edge [bend left=30] node[above,sloped] {\footnotesize$\lambda'_C$} (m-3-3)
(m-2-1) edge [bend right=30] node[below,sloped] {\footnotesize$\lambda'_{C'}$} (m-3-3)
;
\path[dotted,->]
(m-2-2) edge node[description,sloped] {\footnotesize$\exists! \mbox{ } \lambda$} (m-3-3)
;
\end{tikzpicture}
\caption[Universal property of coends]{Universal property of coends.}
\end{figure}
\end{itemize} 
\end{definition}

A coend of a bifunctor $H \colon \C\op \times \C \too \E$, in case it exists, is unique up to isomorphisms, and will be denoted by 
\[
\displaystyle\operatorname*{\int}^{C \in \C} H(C,C).\footnote{The notation $\int H$ will be used for in-text expressions.}
\]
The following result gives us an alternative way to define the coend of a bifunctor, by using coproducts and coequalizers.

\begin{proposition}[equivalent definition of coend]\label{prop:MacLane_alternativo}
Let $H \colon \C\op \times \C \too \E$ a bifunctor, where $\C$ is a small category and $\E$ is a cocomplete category\footnote{Or equivalently, $\E$ has coproducts and coequalizers.}. Then, the coend $\int H$ of $H$ exists, and is given by the coequalizer of the morphisms: 
\begin{align*}
\displaystyle\operatorname*{\coprod}_{f \in \mathsf{Mor}(\C)} H(C',f) \colon \displaystyle\operatorname*{\coprod}_{f \in \mathsf{Mor}(\C)} H(C',C) & \to \displaystyle\operatorname*{\coprod}_{C \in \C} H(C,C), \mbox{ and} \\
\displaystyle\operatorname*{\coprod}_{f \in \mathsf{Mor}(\C)} H(f,C) \colon \displaystyle\operatorname*{\coprod}_{f \in \mathsf{Mor}(\C)} H(C',C) & \to \displaystyle\operatorname*{\coprod}_{C \in \C} H(C,C).
\end{align*}
\end{proposition}

\begin{myproof}
For simplicity, let us write $\alpha = \coprod_{f \in \mathsf{Mor}(\C)} H(C',f)$ and $\beta = \coprod_{f \in \mathsf{Mor}(\C)} H(f,C)$. On the other hand, let ${\rm coeq}(\alpha,\beta)$ denote the coequalizer of $\alpha$ y $\beta$, which is accompanied by a morphism $\lambda \colon \coprod_{C \in \C} H(C,C) \to {\rm coeq}(\alpha,\beta)$ which is natural (that is, $\lambda \circ \alpha = \lambda \circ \beta$) and satisfies the universal property of coequalizers. Let us see that ${\rm coeq}(\alpha,\beta)$ and $H$ satisfy the conditions in \emph{Definition~\ref{def:cofin}}. We first construct morphisms $\lambda_C \colon H(C,C) \to {\rm coeq}(\alpha,\beta)$, for which we verify later the naturality condition.

In what follows, it will come handy to notice that $\alpha$ and $\beta$ are the only morphisms which make the following diagram commute, for the coproduct
\[
\left\{ \displaystyle\operatorname*{\coprod}_{f \in \Hom_{\C}(C,C')} H(C',C), \mbox{ } k_f : H(C',C) \to \displaystyle\operatorname*{\coprod}_{f \in \Hom_{\C}(C,C')} H(C',C) \right\}_{f \in \mathsf{Mor}(\C)}:
\]
\begin{equation}\label{fig:universal_ab} 
\parbox{2.75in}{
\begin{tikzpicture}[description/.style={fill=white,inner sep=2pt}]
\matrix (m) [ampersand replacement=\&, matrix of math nodes, row sep=4em, column sep=3em, text height=2ex, text depth=1.25ex]
{ 
H(C',C') \& \displaystyle\operatorname*{\coprod}_{C \in \C} H(C,C) \\
H(C',C) \& \displaystyle\operatorname*{\coprod}_{f \in \mathsf{Mor}(\C)} H(C',C) \\
H(C,C) \& \displaystyle\operatorname*{\coprod}_{C \in \C} H(C,C) \\
};
\path[->]
($(m-2-1.south)+(0,0.75)$) edge node[left] {\footnotesize$H(C',f)$} ($(m-1-1.north)+(0,-0.75)$)
($(m-2-1.north)+(0,-0.75)$) edge node[left] {\footnotesize$H(f,C)$} ($(m-3-1.south)+(0,0.75)$)
($(m-1-1.east)+(0,0)$) edge node[above] {\footnotesize$j_{C'}$} ($(m-1-2.west)+(0.25,0)$)
($(m-2-1.east)+(0,0)$) edge node[above] {\footnotesize$k_f$} ($(m-2-2.west)+(0.25,0)$)
($(m-3-1.east)+(0,0)$) edge node[above] {\footnotesize$j_C$} ($(m-3-2.west)+(0.25,0)$)
;
\path[dotted,->]
($(m-2-2.south)+(0,0.75)$) edge node[left] {\footnotesize$\exists!$} node[right] {\footnotesize$\alpha$} ($(m-1-2.north)+(0,-0.75)$)
($(m-2-2.north)+(0,-0.75)$) edge node[left] {\footnotesize$\exists!$} node[right] {\footnotesize$\beta$} ($(m-3-2.south)+(0,0.75)$)
;
\end{tikzpicture}
}
\end{equation}
From the natural inclusions $j_C \colon H(C,C) \to \coprod_{C \in \C} H(C,C)$, define:
\[
\lambda_C := \lambda \circ j_C, \mbox{ for every $C \in \C$}.
\]

\newpage

We verify the naturality of the family of morphisms $\{ \lambda_C \}_{C \in \C}$. For each $f \in \Hom_{\C}(C,C')$, we have:
\begin{align*}
\lambda_{C'} \circ H(C',f) & = \lambda \circ j_{C'} \circ H(C',f) = \lambda \circ \alpha \circ k_f = \lambda \circ \beta \circ k_f = \lambda \circ j_C \circ H(f,C) \\
& = \lambda_C \circ H(f,C).
\end{align*}
Now, to show the universal property, suppose we are given another family of morphisms $\lambda'_C \colon H(C,C) \to E$ in $\E$, with $C \in \C$, such that $\lambda'_{C'} \circ H(C',f) = \lambda'_C \circ H(f,C)$ for every $f \in \Hom_{\C}(C,C')$. We find a unique morphism $\omega \colon {\rm coeq}(\alpha,\beta) \to E$ such that $\omega \circ \lambda_C = \lambda'_C$ for every $C \in \C'$. 

First, note that there exists a unique morphism $\lambda' \colon \coprod_{C \in \C} H(C,C) \to E$ such that the following diagram commutes:
\begin{equation}\label{fig:universal_coprod} 
\parbox{2.75in}{
\begin{tikzpicture}[description/.style={fill=white,inner sep=2pt}]
\matrix (m) [ampersand replacement=\&, matrix of math nodes, row sep=3em, column sep=3em, text height=2ex, text depth=1.25ex]
{ 
H(C,C) \& \displaystyle\operatorname*{\coprod}_{C \in \C} H(C,C) \\
{} \& E \\
};
\path[->]
(m-1-1) edge node[above] {\footnotesize$j_C$} (m-1-2)
(m-1-1) edge node[below,sloped] {\footnotesize$\lambda'_C$} (m-2-2)
;
\path[dotted,->]
(m-1-2) edge node[left] {\footnotesize$\exists!$} node[right] {\footnotesize$\lambda'$} (m-2-2)
;
\end{tikzpicture}
}
\end{equation}
The next step is to show that $\lambda' \circ \alpha = \lambda' \circ \beta$. For, we use the universal property of the coproduct $\coprod_{f \in \Mor(\C)} H(C',C)$ in the following commutative diagram:
\[
\begin{tikzpicture}[description/.style={fill=white,inner sep=2pt}]
\matrix (m) [matrix of math nodes, row sep=5em, column sep=5em, text height=1.25ex, text depth=0.25ex]
{ 
  H(C',C) & \displaystyle\operatorname*{\coprod}_{f \in \Mor(\C)} H(C',C) \\
  {} & E \\
};
\path[->]
(m-1-1) edge node[above] {\footnotesize$k_f$} (m-1-2)
(m-1-1) edge node[below,sloped] {\footnotesize$\lambda' \circ j_C \circ H(f,C)$} (m-2-2)
($(m-1-2.north)+(0.1,-1)$) edge node[right] {\footnotesize$\lambda' \circ \beta$} ($(m-2-2.south)+(0.1,0.65)$)
($(m-1-2.north)+(-0.1,-1)$) edge node[left] {\footnotesize$\lambda' \circ \alpha$} ($(m-2-2.south)+(-0.1,0.65)$)
;
\end{tikzpicture}
\]
Indeed, we have
\begin{align*}
(\lambda' \circ \beta) \circ k_f & = \lambda' \circ j_C \circ H(f,C) = \lambda'_C \circ H(f,C) \\
& = \lambda'_{C'} \circ H(C',f) & \mbox{(since the morphisms $\lambda'_C$ are natural)} \\
& = \lambda' \circ j_{C'} \circ H(C',f) \\
& = (\lambda' \circ \alpha) \circ k_f,
\end{align*}
and so $\lambda' \circ \alpha = \lambda' \circ \beta$. Thus, using the universal property of coequalizers, there exists a unique morphism $\omega \colon {\rm coeq}(\alpha,\beta) \to E$ such that the following diagram commutes:
\begin{equation}\label{fig:universal_coeq} 
\parbox{3.75in}{
\begin{tikzpicture}[description/.style={fill=white,inner sep=2pt}]
\matrix (m) [ampersand replacement=\&, matrix of math nodes, row sep=2.5em, column sep=3em, text height=1.5ex, text depth=0.5ex]
{ 
\displaystyle\operatorname*{\coprod}_{f \in \mathsf{Mor}(\C)} H(C',C) \& \displaystyle\operatorname*{\coprod}_{C \in \C} H(C,C) \& {\rm coeq}(\alpha,\beta) \\
{} \& {} \& E \\
};
\path[->]
($(m-1-1.east)+(0,0.1)$) edge node[above,sloped] {\footnotesize$\alpha$} ($(m-1-2.west)+(0,0.1)$)
($(m-1-1.east)+(0,-0.1)$) edge node[below,sloped] {\footnotesize$\beta$} ($(m-1-2.west)+(0,-0.1)$)
(m-1-2) edge node[above] {\footnotesize$\lambda$} (m-1-3)
(m-1-2) edge node[below,sloped] {\footnotesize$\lambda'$} (m-2-3)
; 
\path[dotted,->]
(m-1-3) edge node[left] {\footnotesize$\exists!$} node[right] {\footnotesize$\omega$} (m-2-3)
;
\end{tikzpicture}
}
\end{equation}
From this diagram, the definition of $\lambda_C$, and \eqref{fig:universal_coprod}, it is immediate that $\omega \circ \lambda_C = \lambda'_C$ for every $C \in \C$. Finally, to show that such $\omega$ is unique, suppose there exists $\omega' \colon {\rm coeq}(\alpha,\beta) \to E$ such that $\omega' \circ \lambda_C = \lambda'_C$, for every $C \in \C$. Then, it is easy to see that $(\omega' \circ \lambda) \circ j_C = \lambda' \circ j_C$. By the universal property in \eqref{fig:universal_coprod}, the equality $\omega' \circ \lambda = \lambda'$ holds, and so by the universal property in \eqref{fig:universal_coeq}, we have $\omega' = \omega$.
\end{myproof}

We conclude this section proving the following property of coends.

\begin{proposition}[preservation of coproducts by coends]\label{prop:cofinal_preservacion}
If $\{ H_i \colon \C\op \times \C \too \E \}_{i \in I}$ is a family of functors, where $\E$ is a category with coproducts\footnote{And so coproducts in $\Fun(\C\op \times \C,\E)$ also exist.}, then for the functor $H := \coprod_{i \in I} H_i \colon \C\op \times \C \too \E$, the coend $\int H$ exists and is given by:
\[
\displaystyle\operatorname*{\int}^{C \in \C} H(C,C) = \coprod_{i \in I} \displaystyle\operatorname*{\int}^{C \in \C} H_i(C,C). 
\] 
\end{proposition}

\begin{myproof}
First, it is important to notice, by the universal property of coproducts, that we have the following commutative diagrams:
\begin{equation}\label{fig:Hcovariante} 
\parbox{2.75in}{
\begin{tikzpicture}[description/.style={fill=white,inner sep=2pt}]
\matrix (m) [ampersand replacement=\&, matrix of math nodes, row sep=3em, column sep=4em, text height=1.25ex, text depth=0.25ex]
{ 
  H_i(C',C) \& H(C',C) \\
  H_i(C',C') \& H(C',C') \\
};
\path[->]
(m-1-1) edge node[left] {\footnotesize$H_i(C',f)$} (m-2-1)
;
\path[right hook->]
(m-1-1) edge node[above] {\footnotesize$j^{C',C}_i$} (m-1-2)
(m-2-1) edge node[below] {\footnotesize$j^{C',C'}_i$} (m-2-2)
;
\path[dotted,->]
(m-1-2) edge node[left] {\footnotesize$\exists!$} node[right] {\footnotesize$H(C',f)$} (m-2-2)
;
\end{tikzpicture}
}
\end{equation}

\begin{equation}\label{fig:Hcontravariante} 
\parbox{2.75in}{
\begin{tikzpicture}[description/.style={fill=white,inner sep=2pt}]
\matrix (m) [ampersand replacement=\&, matrix of math nodes, row sep=3em, column sep=4em, text height=1.25ex, text depth=0.25ex]
{ 
  H_i(C',C) \& H(C',C) \\
  H_i(C,C) \& H(C,C) \\
};
\path[->]
(m-1-1) edge node[left] {\footnotesize$H_i(f,C)$} (m-2-1)
;
\path[right hook->]
(m-1-1) edge node[above] {\footnotesize$j^{C',C}_i$} (m-1-2)
(m-2-1) edge node[below] {\footnotesize$j^{C,C}_i$} (m-2-2)
;
\path[dotted,->]
(m-1-2) edge node[left] {\footnotesize$\exists!$} node[right] {\footnotesize$H(f,C)$} (m-2-2)
;
\end{tikzpicture}
}
\end{equation}

\begin{equation}\label{fig:cofinal_de_H} 
\parbox{3.25in}{
\begin{tikzpicture}[description/.style={fill=white,inner sep=2pt}]
\matrix (m) [ampersand replacement=\&, matrix of math nodes, row sep=5em, column sep=3em, text height=0.25ex, text depth=0.25ex]
{ 
  H_i(C,C) \& H(C,C) \\
  \displaystyle\operatorname*{\int}^{C \in \C} H_i(C,C) \& \displaystyle\operatorname*{\coprod}_{i \in I} \displaystyle\operatorname*{\int}^{C \in \C} H_i(C,C) \\
};
\path[->]
($(m-1-1.north)+(0,-0.35)$) edge node[left] {\footnotesize$\lambda^i_C$} ($(m-2-1.south)+(0,0.65)$)
;
\path[right hook->]
(m-1-1) edge node[above] {\footnotesize$j^{C,C}_i$} (m-1-2)
(m-2-1) edge node[below] {\footnotesize$j_i$} (m-2-2)
;
\path[dotted,->]
($(m-1-2.north)+(0.25,-0.35)$) edge node[left] {\footnotesize$\exists!$} node[right] {\footnotesize$\lambda_C$} ($(m-2-2.south)+(0.25,0.65)$)
;
\end{tikzpicture}
}
\end{equation}
where \eqref{fig:Hcovariante} and \eqref{fig:Hcontravariante} show how $H$ is defined on morphisms of $\C\op \times \C$.

\newpage

Having already the object $\coprod_{i \in I} \int H_i$, since $\E$ has coproducts, along with the morphisms $\lambda_C \colon H(C,C) \to \coprod_{i \in I} \int H_i$, we will see that $\coprod_{i \in I} \int H_i$ satisfies the properties in \emph{Definition~\ref{def:coend}}. 
\begin{itemize}[itemsep=2pt,topsep=0pt]
\item[$\bullet$] Let us check that the following diagram commutes (naturality):
\[
\begin{tikzpicture}[description/.style={fill=white,inner sep=2pt}]
\matrix (m) [matrix of math nodes, row sep=3em, column sep=4em, text height=1.25ex, text depth=0.25ex]
{ 
  H(C',C) & H(C',C') \\
  H(C,C) & \displaystyle\operatorname*{\coprod}_{i \in I} \displaystyle\operatorname*{\int}^{C \in \C} H_i(C,C) \\
};
\path[->]
(m-1-1) edge node[above] {\footnotesize$H(C',f)$} (m-1-2)
(m-2-1) edge node[below] {\footnotesize$\lambda_C$} (m-2-2)
(m-1-1) edge node[left] {\footnotesize$H(f,C)$} (m-2-1)
($(m-1-2.north)+(0.25,-0.45)$) edge node[right] {\footnotesize$\lambda_{C'}$} ($(m-2-2.south)+(0.25,0.65)$)
;
\end{tikzpicture}
\]
For this, it suffices to apply the universal property of coproducts in the following diagram: 
\[
\begin{tikzpicture}[description/.style={fill=white,inner sep=2pt}]
\matrix (m) [matrix of math nodes, row sep=5em, column sep=5em, text height=1.25ex, text depth=0.25ex]
{ 
  H_i(C',C) & H(C',C) \\
  {} & \displaystyle\operatorname*{\coprod}_{i \in I} \displaystyle\operatorname*{\int}^{C \in \C} H_i(C,C) \\
};
\path[->]
(m-1-1) edge node[above] {\footnotesize$j^{C',C}_i$} (m-1-2)
($(m-1-1.north)+(0,-0.5)$) edge node[below,sloped] {\footnotesize$\underbrace{j_i \circ \lambda^i_{C'} \circ H_i(C',f)}_{= j_i \circ \lambda^i_C \circ H_i(f,C)}$} ($(m-2-2.south)+(-1.75,0.25)$)
($(m-1-2.north)+(0.1,-0.45)$) edge node[right] {\footnotesize$\lambda_{C} \circ H(f,C)$} ($(m-2-2.south)+(0.1,1.15)$)
($(m-1-2.north)+(-0.1,-0.45)$) edge node[left] {\footnotesize$\lambda_{C'} \circ H(C',f)$} ($(m-2-2.south)+(-0.1,1.15)$)
;
\end{tikzpicture}
\]
Indeed, we have:
\begin{align*}
(\lambda_{C'} \circ H(C',f)) \circ j^{C',C} & = \lambda_{C'} \circ j^{C',C'}_i \circ H_i(C',f) & \mbox{(by~\eqref{fig:Hcovariante})} \\
& = j_i \circ \lambda^i_{C'} \circ H_i(C',f) = j_i \circ \lambda^i_{C} \circ H_i(C,f), & \mbox{(by~\eqref{fig:cofinal_de_H})} \\
(\lambda_C \circ H(f,C)) \circ j^{C',C}_i & = \lambda_C \circ j^{C,C}_i \circ H_i(f,C) & \mbox{(by~\eqref{fig:Hcontravariante})} \\
& = j_i \circ \lambda^i_C \circ H_i(f,C). & \mbox{(by~\eqref{fig:cofinal_de_H})}
\end{align*}

\item[$\bullet$] To check the first part of the universal property of coends, suppose we are given another family $\{ \lambda'_C \colon H(C,C) \to E \}_{C \in \C}$ of natural morphisms in $\E$. We will find a morphism $\omega \colon \coprod_{i \in I} \int H_i \to E$ in $\E$ such that the following diagram commutes:
\[
\begin{tikzpicture}[description/.style={fill=white,inner sep=2pt}]
\matrix (m) [matrix of math nodes, row sep=4em, column sep=3em, text height=1.5ex, text depth=0.25ex]
{ 
  H(C',C) & H(C',C') & {} \\
  H(C,C) & \displaystyle\operatorname*{\coprod}_{i \in I} \displaystyle\operatorname*{\int}^{C \in \C} H_i(C,C) & {} \\
  {} & {} & E \\
};
\path[->]
(m-1-1) edge node[above] {\footnotesize$H(C',f)$} (m-1-2)
(m-2-1) edge node[below] {\footnotesize$\lambda_C$} (m-2-2)
(m-1-1) edge node[left] {\footnotesize$H(f,C)$} (m-2-1)
($(m-1-2.north)+(0.95,-0.25)$) edge [bend left=30] node[above,sloped] {\footnotesize$\lambda'_{C'}$} ($(m-3-3.south)+(0,0.65)$)
($(m-2-1.north)+(0,-0.5)$) edge [bend right=30] node[below,sloped] {\footnotesize$\lambda'_{C}$} ($(m-3-3.south)+(-0.25,0.35)$)
($(m-1-2.north)+(0.25,-0.5)$) edge node[right] {\footnotesize$\lambda_{C'}$} ($(m-2-2.south)+(0.25,0.65)$)
;
\path[dotted,->]
($(m-2-2.north)+(0.5,-0.75)$) edge node[description,sloped] {\footnotesize$\exists! \mbox{ } \omega$} ($(m-3-3.south)+(-0.25,0.5)$)
;
\end{tikzpicture}
\]
On the one hand, for each $i \in I$, we have that:
\begin{align*}
(\eta'_{C'} \circ j^{C',C'}_i) \circ H_i(C',f) & = \eta'_{C'} \circ H(C',f) \circ j^{C',C}_i = \eta'_C \circ H(f,C) \circ j^{C',C}_i & \mbox{(by~\eqref{fig:Hcovariante})} \\
& = (\eta'_C \circ j^{C,C}_i) \circ H_i(f,C). & \mbox{(by~\eqref{fig:Hcontravariante})} 
\end{align*}
Then, for each $i \in I$ there exists a unique morphism $\omega_i \colon \int H_i \to E$ in $\E$ such that the following diagram commutes:
\begin{equation}\label{fig:el_omega_i} 
\parbox{3.25in}{
\begin{tikzpicture}[description/.style={fill=white,inner sep=2pt}]
\matrix (m) [ampersand replacement=\&, matrix of math nodes, row sep=3em, column sep=3em, text height=1.5ex, text depth=0.25ex]
{ 
  H_i(C',C) \& H_i(C',C') \& {} \\
  H_i(C,C) \&\displaystyle\operatorname*{\int}^{C \in \C} H_i(C,C) \& {} \\
  {} \& {} \& E \\
};
\path[->]
(m-1-1) edge node[above] {\footnotesize$H_i(C',f)$} (m-1-2)
(m-2-1) edge node[below] {\footnotesize$\lambda^i_C$} (m-2-2)
(m-1-1) edge node[left] {\footnotesize$H_i(f,C)$} (m-2-1)
($(m-1-2.north)+(0.95,-0.25)$) edge [bend left=30] node[above,sloped] {\footnotesize$\lambda'_{C'} \circ j^{C',C'}_i$} ($(m-3-3.south)+(0,0.65)$)
($(m-2-1.north)+(0,-0.5)$) edge [bend right=30] node[below,sloped] {\footnotesize$\lambda'_{C} \circ j^{C,C}_i$} ($(m-3-3.south)+(-0.25,0.35)$)
($(m-1-2.north)+(0.25,-0.5)$) edge node[right] {\footnotesize$\lambda^i_{C'}$} ($(m-2-2.south)+(0.25,0.65)$)
;
\path[dotted,->]
($(m-2-2.north)+(0.5,-0.75)$) edge node[description,sloped] {\footnotesize$\exists! \mbox{ } \omega_i$} ($(m-3-3.south)+(-0.25,0.5)$)
;
\end{tikzpicture}
}
\end{equation}
Now by the universal property of coproducts we can obtain a unique morphism $\omega \colon \coprod_{i \in I} \int H_i \to E$ in $\E$ making the following diagram commute:
\begin{equation}\label{fig:el_omega} 
\parbox{3in}{
\begin{tikzpicture}[description/.style={fill=white,inner sep=2pt}]
\matrix (m) [ampersand replacement=\&, matrix of math nodes, row sep=3em, column sep=3em, text height=1.5ex, text depth=0.25ex]
{ 
  \displaystyle\operatorname*{\int}^{C \in \C} H_i(C,C) \& \displaystyle\operatorname*{\coprod}_{i \in I} \displaystyle\operatorname*{\int}^{C \in \C} H_i(C,C) \\
  {} \& E \\
};
\path[->]
(m-1-1) edge node[below,sloped] {\footnotesize$\omega_i$} (m-2-2)
;
\path[right hook->]
(m-1-1) edge node[above] {\footnotesize$j_i$} (m-1-2)
;
\path[dotted,->]
(m-1-2) edge node[left] {\footnotesize$\exists!$} node[right] {\footnotesize$\omega$} (m-2-2)
;
\end{tikzpicture}
}
\end{equation}
The equality $\omega \circ \lambda_C = \lambda'_C$, with $C \in \C$, will be a consequence of the commutativity of the diagram below, and the universal property of coproducts:
\[
\begin{tikzpicture}[description/.style={fill=white,inner sep=2pt}]
\matrix (m) [matrix of math nodes, row sep=3em, column sep=4em, text height=1.5ex, text depth=0.25ex]
{ 
  H_i(C,C) & H(C,C) \\
  {} & E \\
};
\path[->]
(m-1-1) edge node[below,sloped] {\footnotesize$\lambda'_C \circ j^{C,C}_i$} (m-2-2)
($(m-1-2.north)+(-0.1,-0.5)$) edge node[left] {\footnotesize$\lambda'_C$} ($(m-2-2.south)+(-0.1,0.65)$)
($(m-1-2.north)+(0.1,-0.5)$) edge node[right] {\footnotesize$\omega \circ \lambda_C$} ($(m-2-2.south)+(0.1,0.65)$)
;
\path[right hook->]
(m-1-1) edge node[above] {\footnotesize$j^{C,C}_i$} (m-1-2)
;
\end{tikzpicture}
\]
Indeed, 
\begin{align*}
& (\omega \circ \lambda_C) \circ j^{C,C}_i = \omega \circ j_i \circ \lambda^i_C = \omega_i \circ \lambda^i_C = \lambda'_C \circ j^{C,C}_i. & \mbox{(by \eqref{fig:cofinal_de_H}, \eqref{fig:el_omega} and \eqref{fig:el_omega_i})}
\end{align*}

\item[$\bullet$] We now prove the uniqueness of $\omega$. Suppose there is another $\omega' \colon  \coprod_{i \in I} \int H_i \to E$ in $\E$ such that $\omega' \circ \lambda_C = \lambda'_C$ for every $C \in \C$. In order to show $\omega' = \omega$, it suffices to verify $\omega' \circ j_i = \omega_i$ and to use the universal property of \eqref{fig:el_omega}. Indeed, we have:
\begin{align*}
(\omega' \circ j_i) \circ \lambda^i_C & = \omega' \circ \lambda_C \circ j^{C,C}_i = \lambda'_C \circ j^{C,C}_i = \omega_i \circ \lambda^i_C. & \mbox{(by \eqref{fig:cofinal_de_H} and \eqref{fig:el_omega_i})} 
\end{align*} 
Therefore, $\int H$ exists and is given by the coproduct $\coprod_{i \in I} \int H_i$. 
\end{itemize}
\end{myproof}


\subsection*{The Kan bifunctor}

In \emph{Theorem~\ref{theo:MacLane}}, for each $D \in \D$ we have a bifunctor 
\[
H_D := \Hom_{\D}(K(-),D) \odot F(-) \colon \C\op \times \C \too \E,
\]
which will be referred to as the \textbf{Kan bifunctor at $\bm{D}$}. On the other hand, since $\E$ is cocomplete, we have that 
\[
\displaystyle\operatorname*{\int}^{C \in \C} H_D(C,C) = \displaystyle\operatorname*{\int}^{C \in \C} \Hom_{\D}(K(C),D) \odot F(C)
\] 
exists by \emph{Proposition~\ref{prop:MacLane_alternativo}}, equipped with natural morphisms $\lambda^D_C \colon H_D(C,C) \to \int H_D$. 

Now suppose that we are given $C, C' \in \C$ and $g \in \Hom_{\D}(D,D')$. Then, we have a function 
\[
\Hom_{\D}(K(C'),g) \colon \Hom_{\D}(K(C'),D) \to \Hom_{\D}(K(C'),D').
\] 
By the universal property in \eqref{fig:t_tensor_E}, there exists a unique morphism 
\[
H_g(C',C) := \Hom_{\D}(K(C'),g) \odot F(C)
\] 
such that the diagram
\begin{equation}\label{fig:bifuntor} 
\parbox{2.75in}{
\begin{tikzpicture}[description/.style={fill=white,inner sep=2pt}]
\matrix (m) [ampersand replacement=\&, matrix of math nodes, row sep=3em, column sep=3em, text height=1.75ex, text depth=0.25ex]
{ 
  F(C) \& \Hom_{\D}(K(C'),D) \odot F(C) \\
  {} \& \Hom_{\D}(K(C'),D') \odot F(C) \\
};
\path[->]
(m-1-1) edge node[above] {\footnotesize$i^{F(C)}_h$} (m-1-2)
(m-1-1) edge node[below,sloped] {\footnotesize$i^{F(C)}_{g \circ h}$} (m-2-2)
;
\path[dotted,->]
(m-1-2) edge node[left] {\footnotesize$\exists!$} node[right] {\footnotesize$H_g(C',C)$} (m-2-2)
;
\end{tikzpicture}
}
\end{equation}
commutes for every $h \in \Hom_{\D}(K(C'),D)$. The morphisms $H_g(C',C)$ satisfy the following properties, which are a straightforward consequence of the universal property of coproducts.

\begin{lemma}\label{lem:propiedad_Hg}
For each pair of morphisms $f \in \Hom_{\C}(C,C')$ and $g \in \Hom_{\D}(D,D')$, the following equalities hold: 
\begin{align*}
H_g(C',C') \circ H_D(C',f) & = H_{D'}(C',f) \circ H_g(C',C), \\
H_g(C,C) \circ H_{D}(f,C) & = H_{D'}(f,C) \circ H_g(C',C). 
\end{align*}
\end{lemma}


\subsection*{The coend of the Kan bifunctor}

Keeping in mind the properties of $H_D$ and $H_g$, consider the mapping $D \mapsto \int H_
D$, with $D \in \D$.

\begin{proposition}\label{prop:Mac_Lane}
The mapping $D \mapsto \int H_D$ gives rise to a functor $\mathfrak{G} \colon \D \too \E$.
\end{proposition}

\begin{myproof}
Define $\mathfrak{G} \colon \D \too \E$ on objects and morphisms of $\D$ as follows:
\begin{itemize}[itemsep=2pt,topsep=0pt]
\item[$\mathsf{(a)}$] $\mathfrak{G}(D) := \int H_D$, for every $D \in \D$.

\item[$\mathsf{(b)}$] Let $g \in \Hom_{\D}(D, D')$. We use the universal property of coends in the following diagram:
\begin{equation}\label{fig:funtor_de_Mac_Lane} 
\parbox{5in}{
\begin{tikzpicture}[description/.style={fill=white,inner sep=2pt}]
\matrix (m) [ampersand replacement=\&, matrix of math nodes, row sep=3em, column sep=5em, text height=1.25ex, text depth=1.25ex]
{ 
  H_D(C',C) \& H_D(C',C') \& H_{D'}(C',C') \\
  H_D(C,C) \& \displaystyle\operatorname*{\int}^{C \in \C} H_D(C,C) \\
  H_{D'}(C,C) \& {} \& \displaystyle\operatorname*{\int}^{C \in \C} H_{D'}(C,C) \\
};
\path[->]
(m-1-1) edge node[above] {\footnotesize$H_D(C',f)$} (m-1-2)
(m-2-1) edge node[below] {\footnotesize$\lambda^{D}_C$} (m-2-2)
(m-1-1) edge node[left] {\footnotesize$H_D(f,C)$} (m-2-1)
(m-1-2) edge node[right] {\footnotesize$\lambda^{D}_{C'}$} (m-2-2)
(m-1-2) edge node[above] {\footnotesize$H_g(C',C')$} (m-1-3)
(m-2-1) edge node[left] {\footnotesize$H_g(C,C)$} (m-3-1)
(m-1-3) edge node[right] {\footnotesize$\lambda^{D'}_{C'}$} (m-3-3)
(m-3-1) edge node[below] {\footnotesize$\lambda^{D'}_C$} (m-3-3)
;
\path[dotted,->]
(m-2-2) edge node[description,sloped] {\footnotesize$\exists! \mbox{ } \mathfrak{G}(g)$} (m-3-3)
;
\end{tikzpicture}
}
\end{equation}
By \emph{Lemma~\ref{lem:propiedad_Hg}}, we have that the outer square in \eqref{fig:funtor_de_Mac_Lane} commutes. Then, there exists a unique morphism $\mathfrak{G}(g) \colon \mathfrak{G}(D) \to \mathfrak{G}(D')$ such that the resulting inner diagrams commute. 
\end{itemize}
By the universal property of coends, one can show that for each pair of morphisms $g_1 \in \Hom_{\D}(D, D')$ and $g_2 \in \Hom_{\D}(D', D'')$, one has $\mathfrak{G}(g_2 \circ g_1) = \mathfrak{G}(g_2) \circ \mathfrak{G}(g_1)$ and $\mathfrak{G}({\rm id}_D) = {\rm id}_{\mathfrak{G}(D)}$. 
\end{myproof}


\subsection*{Proof of the formula to compute Kan extensions}

For the rest of this section, we will focus on proving that $\mathfrak{G}$ is, indeed, the Kan extension of $F$ along $K$. The first thing to do is to construct a natural transformation $\eta \colon F \Rightarrow \mathfrak{G} \circ K$. For each $X \in \C$, we define a morphism $\eta_X \colon F(X) \to \int H_{K(X)}$, in such a way that for each $f \in \Hom_{\C}(X,Y)$ one has the equality 
\begin{align} \label{eqn:Gro_natural}
\mathfrak{G}(K(f)) \circ \eta_X & = \eta_Y \circ F(f).
\end{align}

Recall that $H_{K(X)}(C,C) = \Hom_{\D}(K(C), K(X)) \odot F(C)$. In particular, we have the natural inclusion 
\[
i^{F(X)}_{{\rm id}_{K(X)}} \colon F(X) \hookrightarrow H_{K(X)}(X,X)
\] 
and the natural morphism 
\[
\lambda^{K(X)}_{X} \colon H_{K(X)}(X,X) \to \displaystyle\operatorname*{\int}^{C \in \C} H_{K(X)}(C,C),
\] 
from which we set:
\[
\eta_X := \lambda^{K(X)}_X \circ i^{F(X)}_{{\rm id}_{K(X)}}.
\]
Let us check that \eqref{eqn:Gro_natural} holds:
\begin{align*}
\mathfrak{G}(K(f)) \circ \eta_X & = \mathfrak{G}(K(f)) \circ \lambda^{K(X)}_X \circ i^{F(X)}_{{\rm id}_{K(X)}} \\
& = \lambda^{K(Y)}_{X} \circ H_{K(f)}(X,X) \circ i^{F(X)}_{{\rm id}_{K(X)}} & \mbox{(by~\eqref{fig:funtor_de_Mac_Lane})} \\
& = \lambda^{K(Y)}_{X} \circ i^{F(X)}_{K(f)} & \mbox{(by \eqref{fig:bifuntor})} \\
& = \lambda^{K(Y)}_{X} \circ H_{K(Y)}(f,X) \circ i^{F(X)}_{{\rm id}_{K(Y)}} & \mbox{(by \eqref{fig:t_tensor_E})} \\
& = \lambda^{K(Y)}_{Y} \circ H_{K(Y)}(Y,F(f)) \circ i^{F(X)}_{{\rm id}_{K(Y)}} & \mbox{(by \eqref{fig:funtor_de_Mac_Lane})} \\
& = \lambda^{K(Y)}_{Y} \circ i^{F(Y)}_{{\rm id}_{K(Y)}} \circ F(f) & \mbox{(by \eqref{fig:S_tensor_f})} \\
& = \eta_Y \circ F(f).
\end{align*}

Therefore, we conclude the following result.

\begin{proposition}
The family of morphisms $\{ \eta_X \colon F(X) \to \int H_{K(X)} \}_{X \in \C}$ in $\E$ defines a natural transformation $\eta \colon F \Rightarrow \mathfrak{G} \circ K$, where $\mathfrak{G} \colon \D \too \E$ is the functor from \emph{Proposition~\ref{prop:Mac_Lane}}.
\end{proposition}

We now focus on proving that the pair $(\mathfrak{G},\eta)$ satisfies the universal property of Kan extensions. Suppose we are given a functor $G' \colon \D \too \E$ along with a natural transformation $\eta' \colon F \Rightarrow G' \circ K$. We construct a natural transformation $\alpha \colon \mathfrak{G} \Rightarrow G$ as in \emph{Definition~\ref{def:ext_Kan}}. Each morphism $\alpha_D \colon \int H_D \to G'(D)$ will be obtained using the universal property of coends. First, note that for every $C \in \C$ and $D \in \D$, there exists a unique morphism 
\[
\omega^D_C \colon \Hom_{\D}(K(C),D) \odot F(C) \to G'(D)
\] 
such that the following diagram commutes:
\begin{equation}\label{fig:omegas} 
\parbox{3in}{
\begin{tikzpicture}[description/.style={fill=white,inner sep=2pt}]
\matrix (m) [ampersand replacement=\&, matrix of math nodes, row sep=4em, column sep=3em, text height=1.25ex, text depth=0.25ex]
{ 
  F(C) \& \Hom_{\D}(K(C),D) \odot F(C) \\
  G' \circ K(C) \& G'(D) \\
};
\path[->]
(m-1-1) edge node[above] {\footnotesize$i^{F(C)}_{h}$} (m-1-2)
(m-2-1) edge node[below] {\footnotesize$G'(h)$} (m-2-2)
(m-1-1) edge node[left] {\footnotesize$\eta'_C$} (m-2-1)
;
\path[dotted,->]
(m-1-2) edge node[left] {\footnotesize$\exists!$} node[right] {\footnotesize$\omega^D_C$} (m-2-2)
;
\end{tikzpicture}
}
\end{equation}
where $h \in \Hom_{\D}(K(C),D)$. Now consider the diagram:

\begin{equation}\label{fig:alphas} 
\parbox{4in}{
\begin{tikzpicture}[description/.style={fill=white,inner sep=2pt}]
\matrix (m) [ampersand replacement=\&, matrix of math nodes, row sep=3em, column sep=5em, text height=1.25ex, text depth=0.25ex]
{ 
  H_D(C',C) \& H_{D}(C',C') \\
  H_D(C,C) \& \mathfrak{G}(D) \\
  {} \& {} \& G'(D) \\
};
\path[->]
(m-1-1) edge node[above] {\footnotesize$H_D(C',f)$} (m-1-2)
(m-2-1) edge node[below] {\footnotesize$\lambda^D_C$} (m-2-2)
(m-1-1) edge node[left] {\footnotesize$H_D(f,C)$} (m-2-1)
(m-1-2) edge node[right] {\footnotesize$\lambda^D_{C'}$} (m-2-2)
(m-1-2) edge [bend left=30] node[above,sloped] {\footnotesize$\omega^D_{C'}$} (m-3-3)
(m-2-1) edge [bend right=30] node[below,sloped] {\footnotesize$\omega^D_C$} (m-3-3)
;
\path[dotted,->]
(m-2-2) edge node[description,sloped] {\footnotesize$\exists! \mbox{ } \alpha_D$} (m-3-3)
;
\end{tikzpicture}
}
\end{equation}
Once we prove the equality $\omega^D_{C'} \circ H_D(C',f) = \omega^D_C \circ H_D(f,C)$, that is, that the outer square in \eqref{fig:alphas} commutes, we will be able to assert the existence of a morphism $\alpha_D$ making the resulting inner triangles in \eqref{fig:alphas} commute. 

The equality $\omega^D_{C'} \circ H_D(C',f) = \omega^D_C \circ H_D(f,C)$ will be a consequence of the universal property of coproducts occurring in the following diagram:
\[
\begin{tikzpicture}[description/.style={fill=white,inner sep=2pt}]
\matrix (m) [matrix of math nodes, row sep=3em, column sep=4em, text height=1.75ex, text depth=0.25ex]
{ 
  F(C) & \Hom_{\D}(K(C'),D) \odot F(C) \\
  G' \circ K(C') & G'(D) \\
};
\path[->]
(m-1-1) edge node[above] {\footnotesize$i^{F(C)}_{h}$} (m-1-2)
(m-2-1) edge node[below] {\footnotesize$G'(h)$} (m-2-2)
(m-1-1) edge node[left] {\footnotesize$\eta'_{C'} \circ F(f)$} (m-2-1)
($(m-1-2.east)+(-2.38,-0.25)$) edge node[left] {\footnotesize$\omega^D_{C'} \circ H_D(C',f)$} ($(m-2-2.east)+(-0.75,0.25)$)
($(m-1-2.east)+(-2.18,-0.25)$) edge node[right] {\footnotesize$\omega^D_C \circ H_D(f,C)$} ($(m-2-2.east)+(-0.55,0.25)$)
;
\end{tikzpicture}
\]
where $h \in \Hom_{\D}(K(C'),D)$. Indeed, we have:
\begin{align*}
(\omega^D_{C'} \circ H_D(C',f)) \circ i^{F(C)}_{h} & = \omega^D_{C'} \circ i^{F(C')}_h \circ F(f) & \mbox{(by~\eqref{fig:S_tensor_f})} \\
& = G'(h) \circ \eta'_{C'} \circ F(f) & \mbox{(by~\eqref{fig:omegas})} \\
& = G'(h) \circ G'K(f) \circ \eta'_C & \mbox{(since $\eta'$ is natural)} \\
& = G'(h \circ K(f)) \circ \eta'_C = \omega^D_C \circ i^{F(C)}_{h \circ K(f)} & \mbox{(by~\eqref{fig:omegas})} \\
& = (\omega^D_C \circ H_D(f,C)) \circ i^{F(C)}_{h} & \mbox{(by~\eqref{fig:t_tensor_E})}.
\end{align*}
Hence, $\omega^D_{C'} \circ H_D(C',f) = \omega^D_C \circ H_D(f,C)$. Then, for each $D \in \D$, there exists a unique $\alpha_D \colon \mathfrak{G}(D) \to G'(D)$ as in \eqref{fig:alphas}, that is:
\begin{align}
\omega^D_C & = \alpha_D \circ \lambda^D_C, \label{eqn:igualdad_i} \\
\omega^D_{C'} & = \alpha_D \circ \lambda^D_{C'}. \label{eqn:igualdad_ii}
\end{align}

To show that the family of morphisms $\{ \alpha_D \colon \mathfrak{G}(D) \to G'(D) \}_{D \in \D}$ defines a natural transformation $\alpha \colon \mathfrak{G} \Rightarrow G'$, we use the universal property of coends in the following diagram: 

\begin{equation}\label{fig:alpha_si_es_natural} 
\parbox{4.75in}{
\begin{tikzpicture}[description/.style={fill=white,inner sep=2pt}]
\matrix (m) [ampersand replacement=\&, matrix of math nodes, row sep=3em, column sep=5em, text height=1.25ex, text depth=1.25ex]
{ 
  H_D(C',C) \& H_D(C',C') \& H_{D'}(C',C') \\
  H_D(C,C) \& \mathfrak{G}(D) \\
  H_{D'}(C,C) \& {} \& G'(D') \\
};
\path[->]
(m-1-1) edge node[above] {\footnotesize$H_D(C',f)$} (m-1-2)
(m-2-1) edge node[below] {\footnotesize$\lambda^{D}_C$} (m-2-2)
(m-1-1) edge node[left] {\footnotesize$H_D(f,C)$} (m-2-1)
(m-1-2) edge node[right] {\footnotesize$\lambda^{D}_{C'}$} (m-2-2)
(m-1-2) edge node[above] {\footnotesize$H_g(C',C')$} (m-1-3)
(m-2-1) edge node[left] {\footnotesize$H_g(C,C)$} (m-3-1)
(m-1-3) edge node[right] {\footnotesize$\omega^{D'}_{C'}$} (m-3-3)
(m-3-1) edge node[below] {\footnotesize$\omega^{D'}_C$} (m-3-3)
($(m-2-2.west)+(1,-0.25)$) edge node[below,sloped] {\footnotesize$\alpha_{D'} \circ \mathfrak{G}(g)$} ($(m-3-3.east)+(-1.5,0.25)$)
($(m-2-2.west)+(1.5,-0.25)$) edge node[above,sloped] {\footnotesize$\mathfrak{G}(g) \circ \alpha_D$} ($(m-3-3.east)+(-1.25,0.35)$)
;
\end{tikzpicture}
}
\end{equation}
First, note that the outer square commutes:
\begin{align*}
\omega^{D'}_C \circ H_g(C,C) \circ H_D(f,C) & = \omega^{D'}_C \circ H_{D'}(f,C) \circ H_g(C',C) & \mbox{(by \emph{Lemma~\ref{lem:propiedad_Hg}})} \\
& = \omega^{D'}_{C'} \circ H_{D'}(C',f) \circ H_g(C',C) & \mbox{(by~\eqref{fig:alphas})} \\
& = \omega^{D'}_{C'} \circ H_g(C',C') \circ H_D(C',f) & \mbox{(by \emph{Lemma~\ref{lem:propiedad_Hg}})}.
\end{align*}
Now we show $\alpha_{D'} \circ \mathfrak{G}(g) = \mathfrak{G}(g) \circ \alpha_D$. We will need to use the equality
\begin{align}
\mathfrak{G}(g) \circ \omega^D_C & = \omega^{D'}_C \circ H_g(C,C), \label{eqn:igualdad_iii}
\end{align} 
for every $C \in \C$, which results after applying the universal property of coproducts in the following commutative diagram:
\[
\begin{tikzpicture}[description/.style={fill=white,inner sep=2pt}]
\matrix (m) [matrix of math nodes, row sep=3em, column sep=4em, text height=1.75ex, text depth=0.25ex]
{ 
  F(C) & \Hom_{\D}(K(C),D) \odot F(C) \\
  G' \circ K(C) & G'(D') \\
};
\path[->]
(m-1-1) edge node[above] {\footnotesize$i^{F(C)}_{h}$} (m-1-2)
(m-2-1) edge node[below] {\footnotesize$G'(g \circ h)$} (m-2-2)
(m-1-1) edge node[left] {\footnotesize$\eta'_{C}$} (m-2-1)
($(m-1-2.east)+(-2.38,-0.25)$) edge node[left] {\footnotesize$\mathfrak{G}(g) \circ \omega^D_C$} ($(m-2-2.east)+(-0.86,0.25)$)
($(m-1-2.east)+(-2.18,-0.25)$) edge node[right] {\footnotesize$\omega^{D'}_C \circ H_g(C,C)$} ($(m-2-2.east)+(-0.66,0.25)$)
;
\end{tikzpicture}
\]
where $h \in \Hom_{\D}(K(C),D)$, and:
\begin{align*}
(G'(g) \circ \omega^D_C) \circ i^{F(C)}_h & = G'(g) \circ G'(h) \circ \eta'_C = G'(g \circ h) \circ \eta'_C, & \mbox{(by~\eqref{fig:omegas})} \\
(\omega^{D'}_C \circ H_g(C,C)) \circ i^{F(C)}_h & = \omega^{D'}_C \circ i^{F(C)}_{g \circ h} = G'(g \circ h) \circ \eta'_C. & \mbox{(by \eqref{fig:bifuntor} and \eqref{fig:omegas})} 
\end{align*} 
Thus, we have:
\begin{align*}
(\mathfrak{G}(g) \circ \alpha_D) \circ \lambda^D_C & = \mathfrak{G}(g) \circ \omega^D_C = \omega^{D'}_C \circ H_g(C,C), & \mbox{(by \eqref{eqn:igualdad_i} and \eqref{eqn:igualdad_iii})} \\
(\alpha_{D'} \circ \mathfrak{G}(g)) \circ \lambda^D_C & = \alpha_{D'} \circ \lambda^{D'}_C = \omega^{D'}_C \circ H_g(C,C) \circ H_g(C,C) & \mbox{(by~\eqref{fig:funtor_de_Mac_Lane} and \eqref{eqn:igualdad_i})}.
\end{align*}
In a similar way, we have: 
\begin{align*} 
(\mathfrak{G}(g) \circ \alpha_D) \circ \lambda^D_{C'} & = \omega^{D'}_C \circ H_g(C,C), \\ 
(\alpha_{D'} \circ \mathfrak{G}) \circ \lambda^{D}_{C'} & = \omega^{D'}_{C'} \circ H_g(C',C').
\end{align*}
Hence, $\alpha_{D'} \circ \mathfrak{G}(g) = \mathfrak{G}(g) \circ \alpha_D$, that is, $\alpha \colon \mathfrak{G} \Rightarrow G'$ is a natural transformation.   

Now we prove that $\eta' = \alpha_K \circ \eta$, that is, $\eta'_X = \alpha_{K(X)} \circ \eta_X$ for every $X \in \C$:
\begin{align*}
\alpha_{K(X)} \circ \eta_X & = \alpha_{K(X)} \circ \lambda^{K(X)}_X \circ i^{F(X)}_{{\rm id}_{K(X)}} = \omega^{K(X)}_X \circ i^{F(X)}_{{\rm id}_{K(X)}} & \mbox{(by~\eqref{eqn:igualdad_i})} \\
& = \mathfrak{G}({\rm id}_{K(X)}) \circ \eta'_X & \mbox{(by~\eqref{fig:omegas})} \\
& = \eta'_X.
\end{align*}

Finally, let us show the uniqueness of $\alpha$. Suppose there exists another natural transformation $\alpha' \colon \mathfrak{G} \Rightarrow G'$ such that $\eta' = \alpha'_K \circ \eta$. For each $D \in \D$, the equality $\alpha_D = \alpha'_D$ will follow from the universal property in \eqref{fig:alphas}, after verifying that each $\alpha'_D$ satisfies the equalities 
\begin{align}
\alpha'_D \circ \lambda^D_C & = \omega^D_C, \label{eqn:igualdad_prima_i} \\
\alpha'_D \circ \lambda^D_{C'} & = \omega^D_{C'}. \label{eqn:igualdad_prima_ii}
\end{align}
To show \eqref{eqn:igualdad_prima_i}, we check that the following diagram commutes:
\[
\begin{tikzpicture}[description/.style={fill=white,inner sep=2pt}]
\matrix (m) [matrix of math nodes, row sep=3em, column sep=4em, text height=1.75ex, text depth=0.25ex]
{ 
  F(C) & \Hom_{\D}(K(C),D) \odot F(C) \\
  G'K(C) & G'(D) \\
};
\path[->]
(m-1-1) edge node[above] {\footnotesize$i^{F(C)}_{h}$} (m-1-2)
(m-2-1) edge node[below] {\footnotesize$G'(h)$} (m-2-2)
(m-1-1) edge node[left] {\footnotesize$\eta'_{C}$} (m-2-1)
($(m-1-2.east)+(-2.38,-0.25)$) edge node[left] {\footnotesize$\omega^D_C$} ($(m-2-2.east)+(-0.81,0.25)$)
($(m-1-2.east)+(-2.18,-0.25)$) edge node[right] {\footnotesize$\alpha'_D \circ \lambda^D_C$} ($(m-2-2.east)+(-0.61,0.25)$)
;
\end{tikzpicture}
\]
where $h \in \Hom_{\D}(K(C),D)$. We have: 
\begin{align*}
(\alpha'_D \circ \lambda^D_C) \circ i^{F(C)}_h & = \alpha'_D \circ \lambda^D_C \circ H_h(C,C) \circ i^{F(C)}_{{\rm id}_{K(C)}} & \mbox{(by~\eqref{fig:bifuntor})} \\
& = \alpha'_D \circ \mathfrak{G}(h) \circ \lambda^{K(C)}_C \circ i^{F(C)}_{{\rm id}_{K(C)}} & \mbox{(by~\eqref{fig:funtor_de_Mac_Lane})} \\
& = G'(h) \circ \alpha'_{K(C)} \circ \lambda^{K(C)}_C \circ i^{F(C)}_{{\rm id}_{K(C)}} & \mbox{(since $\alpha'$ is a natural transformation)} \\
& = G'(h) \circ \alpha'_{K(C)} \circ \eta_C \\
& = G'(h) \circ \eta'_C.
\end{align*}
Then, \emph{\eqref{eqn:igualdad_prima_i}} holds, while \eqref{eqn:igualdad_prima_ii} follows in a similar way. Thus, we conclude $\alpha = \alpha'$. Therefore, $(\mathfrak{G},\eta)$ is the Kan extension of $F$ along $K$, thus proving \emph{Theorem~\ref{theo:MacLane}}.


\section{Kan extensions and the tensor product of $\C$-modules}\label{sec:producto}

Suppose we are given two functors $F \colon \D\op \too \mathcal{V}$ and $G \colon \D \too \mathcal{
M}$, and a bifunctor $- \otimes - \colon \mathcal{V} \times \mathcal{M} \too \mathcal{M}$. If the coend of the bifunctor $F(-) \otimes G(-) \colon \D\op \times \D \too \mathcal{M}$ exists, the \textbf{tensor product} of $F$ and $G$ is defined as:
\[
F \displaystyle\operatorname*{\otimes}_{\D} G := \displaystyle\operatorname*{\int}^{D \in \D} F(D) \otimes G(D).
\]
The following result, which is a consequence of \emph{Proposici\'on~\ref{prop:MacLane_alternativo}}, establishes some conditions under which the tensor product $F \otimes_{\D} G$ exists, along with an alternative way to define it.

\begin{corollary}[equivalent definition of tensor product of functors]\label{coro:existencia_producto}
Let $\D$ be a small category and $\mathcal{M}$ be cocomplete. Then, the tensor product $F \otimes_{\D} G$ exists, and is given by the coequalizer of the following morphisms:
\begin{align*}
\displaystyle\operatorname*{\coprod_{g \in \mathsf{Mor}(\D)}} F(D') \otimes G(g) \colon \coprod_{g \in \Mor(\D)} F(D') \otimes G(D) & \to \coprod_{D \in \D} F(D) \otimes G(D), \mbox{ and} \\
\displaystyle\operatorname*{\coprod_{g \in \mathsf{Mor}(\D)}} F(g) \otimes G(D) \colon \coprod_{g \in \Mor(\D)} F(D') \otimes G(D) & \to \coprod_{D \in \D} F(D) \otimes G(D).
\end{align*}
\end{corollary}

We motivate this section by saying that, in the case $\mathcal{M} = \mathcal{V} = \Set$, there is another way to compute tensor products of functors via adjunctions. Specifically, in this case it is known that $- \otimes_{\D} -$ has a left adjoint. The reader can find the proof of this fact in \emph{\cite[Chapter VII, Theorem 1]{MML}}. Then, by the dual of \emph{Theorem~\ref{theo:Kan_extensions_vs_adjoint_pairs}}, we have that $- \otimes_{\D} -$ is a Kan extension.

\begin{example}\label{ej:producto_tensorial_funtores} 
We state some examples of tensor products of functors:
\begin{itemize}[itemsep=2pt,topsep=0pt]
\item[$\uno$] Let $R$ be an associative ring with identity. Let $A \in \Mod(R\op)$ and $B \in \Mod(R)$. Note that $A$ and $B$ can be regarded as functors $A \colon \mathfrak{R}\op \too \Ab$ and $B \colon \mathfrak{R} \too \Ab$, where $\mathfrak{R}$ is the category with only one object, and whose morphisms are given by the elements of $R$. On the other hand, we have a bifunctor $- \otimes_{\mathbb{Z}} - \colon \Ab \times \Ab \too \Ab$, given by the standard tensor product of abelian groups. In this case, the tensor product of modules $A \otimes_R B$ is isomorphic to the tensor product of functors $A \otimes_{\mathfrak{R}} B$. 

\item[$\dos$] Consider a pair of functors $F \colon \C \too \E$ and $K \colon \C \too \D$, where $\C$ is a small category and $\E$ is cocomplete. For $D \in \D$ fixed, set the functor
\[
G := \Hom_{\D}(K(-),D) \colon \C^{\rm op} \too \Set.
\] 
In this case, we have:
\begin{align*}
G \displaystyle\operatorname*{\otimes}_{\C} F & = \displaystyle\operatorname*{\int}^{C \in \C} \Hom_{\D}(K(C),D) \odot F(C) & \mbox{(by the definition of tensor product of functors)} \\
& = \Lan_K(F)(D) & \mbox{(by \emph{Theorem~\ref{theo:MacLane}})}.
\end{align*}
In particular, we have the equality:
\begin{align}
\Hom_{\D}(-,D) \displaystyle\operatorname*{\otimes}_{\D} F & = \Lan_{{\rm id}}(F)(D) = F(D). \label{eqn:derecha}
\end{align}
\end{itemize}
\end{example}

We devote the rest of this section to studying the tensor product of $\C$-modules, and to showing that it represents another example of tensor products of functors.


\subsection*{The category of $\C$-modules}

Let $\C$ be a skeletally small and preadditive category. Denote by $\Mod(\C) := \Fun(\C\op, \Ab)$ the category of contravariant functors from $\C$ to $\Ab$. We will refer to $\Mod(\C)$ as the \textbf{category of} (\textbf{right}) \textbf{$\bm{\C}$-modules}. This category is studied, for example, B. Mitchell's book  \emph{\cite[page 106]{Mitchell}}. From this reference, we collect below some properties of $\Mod(\C)$:
\begin{itemize}[itemsep=2pt,topsep=0pt]
\item[$\bullet$] $\Mod(\C)$ is an abelian category, since $\Ab$ is abelian, where the zero object is given by the zero functor $F \colon \C\op \too \Ab$ defined as $F(C) := 0$, for every $C \in \C$. 

\item[$\bullet$] The kernel of a morphism $\alpha \colon F \Rightarrow G$ in $\Mod(\C)$ is defined component-wise, that is, for each $f \in \Hom_{\C}(C_1,C_2)$, $\Ker(\alpha)(f)$ is the only morphism $h \colon \Ker(\alpha_{C_2}) \to \Ker(\alpha_{C_1})$ in $\Ab$ making the following diagram commute:
\[
\begin{tikzpicture}[description/.style={fill=white,inner sep=2pt}]
\matrix (m) [matrix of math nodes, row sep=3em, column sep=4em, text height=1.25ex, text depth=0.25ex]
{ 
  \Ker(\alpha_{C_2}) & F(C_2) & G(C_2) \\
  \Ker(\alpha_{C_1}) & F(C_1) & G(C_1) \\
};
\path[right hook->]
(m-1-1) edge (m-1-2) 
(m-2-1) edge (m-2-2) 
;
\path[->]
(m-1-2) edge node[above] {\footnotesize$\alpha_{C_2}$} (m-1-3)
(m-2-2) edge node[below] {\footnotesize$\alpha_{C_1}$} (m-2-3)
(m-1-2) edge node[right] {\footnotesize$F(f)$} (m-2-2)
(m-1-3) edge node[right] {\footnotesize$G(f)$} (m-2-3)
;
\path[dotted,->]
(m-1-1) edge node[left] {\footnotesize$\exists!$} node[right] {\footnotesize$h$} (m-2-1)
;
\end{tikzpicture}
\]
It is not difficult to show that the left-hand square in the previous diagram gives rise to a natural transformation $\Ker(\alpha) \Rightarrow F$ which satisfies the universal property of kernels. Dually, one can show that the category $\Mod(\C)$ also has cokernels defined component-wise.

\item[$\bullet$] $\Mod(\C)$ has arbitrary coproducts. Suppose we are given a family of $\C$-modules $\{ F_i \}_{i \in I}$. Then, for each $C \in \C$, $\{ F_i(C) \}_{i \in I}$ is a family of abelian groups, for which there exists the coproduct $\coprod_{i \in I} F_i(C)$, along with natural inclusions $j^i_C \colon F_i(C) \to \coprod_{i \in I} F_i(C)$. Now suppose that we have a morphism $f \in \Hom_{\C}(C_1,C_2)$. By the universal property of coproducts in $\Ab$, there exists a unique morphism 
\[
\coprod_{i \in I} F_i(f) \colon \coprod_{i \in I} F_i(C_2) \to \coprod_{i \in I} F_i(C_1)
\] 
in $\Ab$ such that the following diagram commutes:
\[
\begin{tikzpicture}[description/.style={fill=white,inner sep=2pt}]
\matrix (m) [matrix of math nodes, row sep=3em, column sep=4em, text height=1.25ex, text depth=1ex]
{ 
  F_i(C_2) & \displaystyle\operatorname*{\coprod}_{i \in I} F_i(C_2) \\
  F_i(C_1) & \displaystyle\operatorname*{\coprod}_{i \in I} F_i(C_1) \\
};
\path[right hook->]
(m-1-1) edge node[above] {\footnotesize$j^i_{C_2}$} (m-1-2) 
(m-2-1) edge node[below] {\footnotesize$j^i_{C_1}$} (m-2-2) 
;
\path[->]
(m-1-1) edge node[left] {\footnotesize$F_i(f)$} (m-2-1)
;
\path[dotted,->]
(m-1-2) edge node[left] {\footnotesize$\exists!$} node[right] {\footnotesize$\displaystyle\operatorname*{\coprod}_{i \in I} F_i(f)$} (m-2-2)
;
\end{tikzpicture}
\]
The previous diagram is functorial on the morphisms of $\C$, and so we can define the coproduct $\C$-module $\coprod_{i \in I} F_i \colon \C\op \too \Ab$ as follows:
\[
\left( \displaystyle\operatorname*{\coprod}_{i \in I} F_i \right)(f) := \displaystyle\operatorname*{\coprod}_{i \in I} F_i(f).
\]
Moreover, for each $i \in I$, we have a natural transformation $j^i \colon F_i \Rightarrow \coprod_{i \in I} F_i$ defined by $(j^i)_C := j^i_C$. The $\C$-module $\coprod_{i \in I} F_i$, along with the family of natural transformation $j^i \colon F_i \Rightarrow \coprod_{i \in I} F_i$, defines the coproduct of the family of $\C$-modules $\{ F_i \}_{i \in I}$. Dually, one can show that $\Mod(\C)$ has arbitrary products. 

\item[$\bullet$] Every monomorphism in $\Mod(\C)$ is the kernel of some morphism in $\Mod(\C)$. Indeed, observe that $\alpha \colon F \Rightarrow G$ is a monomorphism in $\Mod(\C)$ if, and only if, $\Ker(\alpha) = 0$, which in turn is equivalent to $\Ker(\alpha_C) = 0$, for every $C \in \C$. From this we can deduce that if $\alpha \colon F \Rightarrow G$ is a monomorphism in $\Mod(\C)$, then for each $f \in \Hom_{\C}(C_1,C_2)$ we have the following commutative diagram: 
\[
\begin{tikzpicture}[description/.style={fill=white,inner sep=2pt}]
\matrix (m) [matrix of math nodes, row sep=3.5em, column sep=3em, text height=1.25ex, text depth=0.25ex]
{ 
0 & F(C_2) & G(C_2) & \Coker(\alpha_{C_2}) & 0 \\
0 & F(C_1) & G(C_1) & \Coker(\alpha_{C_1}) & 0 \\
};
\path[->]
(m-1-1) edge (m-1-2) (m-1-2) edge node[above] {\footnotesize$\alpha_{C_2}$} (m-1-3) (m-1-3) edge node[above] {\footnotesize$\pi_{C_2}$} (m-1-4) (m-1-4) edge (m-1-5)
(m-2-1) edge (m-2-2) (m-2-2) edge node[below] {\footnotesize$\alpha_{C_1}$} (m-2-3) (m-2-3) edge node[below] {\footnotesize$\pi_{C_1}$} (m-2-4) (m-2-4) edge (m-2-5)
(m-1-2) edge node[left] {\footnotesize$F(f)$} (m-2-2) 
(m-1-3) edge node[left] {\footnotesize$G(f)$} (m-2-3)
;
\path[dotted,->]
(m-1-4) edge node[left] {\footnotesize$\exists!$} node[right] {\footnotesize$\Coker(\alpha)(f)$} (m-2-4)
;
\end{tikzpicture}
\]
One can see that the family of morphisms $\pi_C \colon G(C) \to \Coker(\alpha_C)$ defines a natural transformation $\pi \colon G \Rightarrow \Coker(\alpha)$, and that $\Ker(\pi) = \alpha$. Dually, one can note that every epimorphism in $\Mod(\C)$ is the cokernel of a morphism in $\Mod(\C)$.

\item[$\bullet$] Since kernels and cokernels of morphisms in $\Mod(\C)$ are defined component-wise, we have that a sequence of $\C$-modules
\[
M_1 \Rightarrow M_2 \Rightarrow M_3
\] 
is exact if, and only if, for each $C \in \C$, the sequence
\[
M_1(C) \to M_2(C) \to M_3(C)
\]
is exact in $\Ab$.
\end{itemize}

Probably the reader has already noted that the category $\Mod(\C)$ has an structure which is richer than that of an abelian category. Besides having arbitrary products and coproducts, $\Mod(\C)$ is also equipped with enough projective and injective objects. 
\begin{itemize}[itemsep=2pt,topsep=0pt]
\item[$\bullet$] Every $\C$-module of the form $\Hom_{\C}(-,C) \colon \C\op \too \Ab$, with $C \in \C$ fixed, is a projective object in $\Mod(\C)$. Moreover, the mapping $C \mapsto \Hom_{\C}(-,C)$ defines a (fully faithful) functor $\C \too \Mod(\C)$, known as the \emph{Yoneda embedding}.  

\item[$\bullet$] For each family $\{ C_i \}_{i \in I}$ of objects of $\C$, one has that the coproduct $\coprod_{i \in I} \Hom_{\C}(-,C_i)$ is a projective $\C$-module. We will refer to this type of $\C$-modules as \textbf{free}. Moreover, $\Mod(\C)$ has enough projective objects, due to the fact that for each $\C$-module $M$, one can always construct an epimorphism of the form $\coprod_{i \in I} \Hom_{\C}(-,C_i) \Rightarrow M$, for some family of objects $\{ C_i \}_{i \in I}$ in $\C$. 
\end{itemize}


\subsection*{Tensor product of $\C$-modules}

The tensor product of $\C$-modules is a fundamental construction in tilting theory, as part of representation theory of algebras. In the 70s, several works by M. Auslander triggered the importance of the study of the category of functors (among them, the category of $\C$-modules) and its tensor product in some contexts of tilting theory of Artin algebras. Among these works one can find, for instance, \emph{\cite{Auslander}}. Recently, R. Mart\'inez-Villa and M. Ortiz-Morales revisit Auslander's works to study tilting theory in the category $\Mod(\C)$ of $\C$-modules, constructing in \emph{\cite{Martin}} a torsion class from a tilting subcategory of $\Mod(\C)$, and also proving properties of such classes relative to the tensor product in $\Mod(\C)$.

Keeping in mind the importance of the tensor product of $\C$-modules, we devote this section to giving a more categorical approach to this notion by using the coend construction. 

We begin recalling how to define the tensor product of $\C$-modules. First, suppose we are given a $\C\op$-module $F \colon \C \too \Ab$. Next, we will construct a functor
\[
F \otimes - \colon \Mod(\C) \too \Ab.
\] 
Let $\A$ be the full subcategory of $\Mod(\C)$ whose objects are the free $\C$-modules. Define $F \otimes - \colon \A \too \Ab$ as follows:
\begin{itemize}[itemsep=2pt,topsep=0pt]
\item[$\iroman$] \underline{Definition of $F \otimes -$ on objects of $\A$}: For each object $\coprod_{i \in I} \Hom_{\C}(-,C_i)$ of $\A$, set
\[
F \otimes \left(\displaystyle\operatorname*{\coprod}_{i \in I} \Hom_{\C}(-,C_i)\right) := \displaystyle\operatorname*{\coprod}_{i \in I} F(C_i).
\] 
By \emph{Yoneda Lemma}, this expression is well defined.

\item[$\iiroman$] \underline{Definiton of $F \otimes -$ on morphisms of $\A$}: For each morphism in $\A$, that is, a morphism of the form 
\[
( \Hom_{\C}(-,f_{ij}) )_{(i,j) \in I \times J} \colon \coprod_{i \in I} \Hom_{\C}(-,C_i) \Rightarrow \coprod_{j \in J} \Hom_{\C}(-,C_j),
\] 
where $f_{ij} \in \Hom_\C(C_i,C_j)$, define
\[
F \otimes \left( \Hom_{\C}(-,f_{ij}) \right)_{(i,j) \in I \times J} := (F(f_{ij}))_{(i,j) \in I \times J},
\] 
as the morphism induced by the universal property of coproducts.  
\end{itemize}
We use partial free resolutions to extend $F \otimes -$ on the whole category $\Mod(\C)$. 
\begin{itemize}[itemsep=2pt,topsep=0pt]
\item[$\mathsf{(iii)}$] \underline{Definition of $F \otimes -$ on objects of $\Mod(\C)$}: Let $G \in \Mod(\C)$, and consider a partial free resolution of $G$ of length 1, that is, an exact sequence in $\Mod(\C)$ of the form
\[
\coprod_{i \in I} \Hom_{\C}(-,C_i) \xRightarrow{(\Hom_{\C}(-,f_{ij}))_{(i,j) \in I \times J}} \coprod_{j \in J} \Hom_{\C}(-,C_j) \Rightarrow G \Rightarrow 0.
\]
Applying the functor $F \otimes -$ to the left-hand morphism, we obtain the following exact sequence in $\Ab$:
\[
\coprod_{i \in I} F(C_i) \xrightarrow{(F(f_{ij}))_{(i,j) \in I \times J}} \coprod_{j \in J} F(C_j) \to \Coker((F(f_{ij}))_{(i,j) \in I \times J}) \to 0.
\]
Thus, define
\[
F \otimes G := \Coker((F(f_{ij}))_{(i,j) \in I \times J}).
\]

\item[$\mathsf{(iv)}$] \underline{Definition of $F \otimes -$ on morphisms of $\Mod(\C)$}: Let $\eta \colon G \Rightarrow H$ be a morphism of $\C$-modules. Using the universal property of coproducts, one can find unique morphisms in $\Mod(\C)$
\[
\eta' \colon \coprod_{j \in J} \Hom_{\C}(-,C_j) \Rightarrow \coprod_{j \in J'} \Hom_{\C}(-,C'_j) \mbox{ \ y \ } \eta'' \colon \coprod_{i \in I} \Hom_{\C}(-,C_i) \Rightarrow \coprod_{i \in I'} \Hom_{\C}(-,C'_i)
\] 
such that the following diagram commutes:
\[
\begin{tikzpicture}[description/.style={fill=white,inner sep=2pt}]
\matrix (m) [matrix of math nodes, row sep=3.5em, column sep=3em, text height=1.25ex, text depth=0.25ex]
{ 
\displaystyle\operatorname*{\coprod}_{i \in I} \Hom_{\C}(-,C_i) & {} & {} & \displaystyle\operatorname*{\coprod}_{j \in J} \Hom_{\C}(-,C_j) & G & 0 \\
\displaystyle\operatorname*{\coprod}_{i \in I'} \Hom_{\C}(-,C'_i) & {} & {} & \displaystyle\operatorname*{\coprod}_{j \in J'} \Hom_{\C}(-,C'_j) & H & 0 \\
};
\path[->]
(m-1-1) edge [double,double equal sign distance,-implies] node[above] {\footnotesize$(\Hom_{\C}(-,f_{ij}))_{(i,j) \in I \times J}$} (m-1-4) (m-1-4) edge [double,double equal sign distance,-implies] (m-1-5) (m-1-5) edge [double,double equal sign distance,-implies] (m-1-6)
(m-2-1) edge [double,double equal sign distance,-implies] node[below] {\footnotesize$(\Hom_{\C}(-,f'_{ij}))_{(i,j) \in I' \times J'}$} (m-2-4) (m-2-4) edge [double,double equal sign distance,-implies] (m-2-5) (m-2-5) edge [double,double equal sign distance,-implies] (m-2-6)
(m-1-5) edge [double,double equal sign distance,-implies] node[left] {\footnotesize$\eta$} (m-2-5)
;
\path[dotted,->]
(m-1-1) edge [double,double equal sign distance,-implies] node[left] {\footnotesize$\eta''$} (m-2-1)
(m-1-4) edge [double,double equal sign distance,-implies] node[left] {\footnotesize$\eta'$} (m-2-4)
;
\end{tikzpicture}
\]
Applying the functor $F \otimes - \colon \A \too \Ab$ to the previous diagram, we obtain the following commutative diagram in $\Ab$:
\[
\begin{tikzpicture}[description/.style={fill=white,inner sep=2pt}]
\matrix (m) [matrix of math nodes, row sep=3.5em, column sep=2em, text height=1.25ex, text depth=0.25ex]
{ 
\displaystyle\operatorname*{\coprod}_{i \in I} F(C_i) & {} & {} & \displaystyle\operatorname*{\coprod}_{j \in J} F(C_j) & \overbrace{\Coker((F(f_{ij}))_{(i,j) \in I \times J})}^{= F \otimes G} & 0 \\
\displaystyle\operatorname*{\coprod}_{i \in I'} F(C'_i) & {} & {} & \displaystyle\operatorname*{\coprod}_{j \in J'} F(C'_j) & \underbrace{\Coker((F(f'_{ij}))_{(i,j) \in I' \times J'})}_{= F \otimes H} & 0 \\
};
\path[->]
(m-1-1) edge node[above] {\footnotesize$(F(f_{ij}))_{(i,j) \in I \times J}$} (m-1-4)
(m-2-1) edge node[below] {\footnotesize$(F(f'_{ij}))_{(i,j) \in I' \times J'}$} (m-2-4)
(m-1-1) edge node[left] {\footnotesize$F \otimes \eta''$} (m-2-1)
(m-1-4) edge node[left] {\footnotesize$F \otimes \eta'$} (m-2-4)
(m-1-4) edge (m-1-5) (m-1-5) edge (m-1-6)
(m-2-4) edge (m-2-5) (m-2-5) edge (m-2-6)
;
\path[dotted,->]
(m-1-5) edge node[left] {\footnotesize$\exists!$} node[right] {\footnotesize$h$} (m-2-5)
;
\end{tikzpicture}
\]
Thus, define $F \otimes \eta \colon F \otimes G \to F \otimes H$ as the only morphism $h$ appearing in the previous diagram. 
\end{itemize}

Dually, given a $\C$-module $G \colon \C\op \too \Ab$, we can construct a functor 
\[
- \otimes G \colon \Mod(\C\op) \too \Ab
\] 
such that 
\[
\left( \coprod_{i \in I} \Hom_{\C}(C_i,-) \right) \otimes G = \coprod_{i \in I} G(C_i),
\]
for every family $\{ C_i \}_{i \in I}$ of objects of $\C$. The properties of the functors $F \otimes -$ and $- \otimes G$ are summarized in the following result, which is known from \emph{\cite{Auslander}}.

\begin{theorem}[tensor product of $\C$-modules]\label{theo:producto_tensorial}
Given a skeletally small preadditive category $\C$, there exists a unique (up to natural isomorphisms) bifunctor 
\[
- \otimes - : \Mod(\C\op) \times \Mod(\C) \too \Ab,
\] 
called \textbf{tensor product}, such that for each $F \in \Mod(\C\op)$ and $G \in \Mod(\C)$, the following properties hold: 
\begin{itemize}[itemsep=2pt,topsep=0pt]
\item[$\uno$] The functors $F \otimes - : \Mod(\C) \too \Ab$ and $- \otimes G : \Mod(\C\op) \too \Ab$ are right exact and preserve arbitrary coproducts in $\Mod(\C)$ and $\Mod(\C\op)$, respectively. 

\item[$\dos$] $F \otimes \Hom_{\C}(-,C) = F(C)$ and $\Hom_{\C}(C,-) \otimes G = G(C)$, for every $C \in \C$.
\end{itemize}
\end{theorem}

Keeping in mind this review of $\C$-modules, it is time to see how the tensor product of $\C$-modules $- \otimes -$ is naturally isomorphic to the tensor product of functors $- \otimes_{\C\op} -$. First, we will show that $- \otimes_{\C\op} -$ satisfies conditions $\uno$ and $\dos$ in \emph{Theorem~\ref{theo:producto_tensorial}}. So $- \otimes_{\C\op} -$ and $- \otimes -$ will be naturally isomorphic. In the end of these notes, we will define an explicit natural isomorphism $- \otimes - \Rightarrow - \otimes_{\C\op} -$.


\subsection*{The tensor product of $\C$-modules as a tensor product of functors}

Consider the standard tensor product of abelian groups $- \otimes_{\mathbb{Z}} - : \Ab \times \Ab \too \Ab$. Let us show that for every $F \in \Mod(\C\op)$ and $G \in \Mod(\C)$, the functor tensor product $F \otimes_{\C\op} G := \int F \otimes_{\mathbb{Z}} G$ is isomorphic to the tensor product $F \otimes G$ from \emph{Theorem~\ref{theo:producto_tensorial}}. Before proving this fact, we need several results. The first thing to notice is that $F \otimes_{\C\op} G$ does not necessarily exist, since according to \emph{Corollary~\ref{coro:existencia_producto}}, the category $\C$ needs to be small. We will see that this requisite can be replaced by the more general setting in whicn $\C$ is skeletally small. The tensor product of functors $F \otimes_{\C\op} G$ can be computed by using the coend formula \eqref{eqn:Grothendieck} if $\C$ is small. In order to extend this computation for $\C$ skeletally small, we prove the following extension of \emph{Theorem~\ref{theo:MacLane}}.

\begin{proposition}\label{prop:MacLaneSke}
Let $F \colon \C \too \E$ and $K \colon \C \too \D$ be two functors, and $J \colon \C' \too \C$ an equivalence of categories such that $\Lan_{K \circ J}(F \circ J)$ exists. Then so does $\Lan_K(F)$, and $\Lan_K(F) = \Lan_{K \circ J}(F \circ J)$. 
\end{proposition}

\begin{myproof}
Let $P \colon \C \too \C'$ be a functor such that there are natural isomorphisms $\mu \colon {\rm id}_{\C'} \Rightarrow P \circ J$ and $\nu \colon J \circ P \Rightarrow {\rm id}_{\C}$. Since $J$ is an equivalence of categories, $(J,P)$ is an adjoint pair, and hence $\mu$ and $\nu$ satisfy the triangle identities \eqref{eqn:triangular_uno} and \eqref{eqn:triangular_dos}. Denote $G = \Lan_{K \circ J}(F \circ J)$, accompanied by a natural transformation $\gamma \colon F \circ J \Rightarrow G \circ K \circ J$ which satisfies the universal property of Kan extensions. 

In what follows, we show that $G = \Lan_K(F)$. We need to obtain a natural transformation $\eta \colon F \Rightarrow G \circ K$ as in \emph{Definition~\ref{def:ext_Kan}}. We define such $\eta$ as:
\[
\eta := G \circ K(\nu) \circ \gamma_P \circ F(\nu^{-1}). 
\]
Suppose we are given another functor $G' \colon \D \too \E$ along with a natural transformation $\eta' \colon F \Rightarrow G' \circ K$. Since $G = \Lan_{K \circ J}(F \circ J)$, there exists a unique $\alpha \colon G \Rightarrow G'$ such that $\eta'_J = \alpha_{K \circ J} \circ \gamma$. Let us see $\eta' = \alpha_K \circ \eta$. Let $C \in \C$. Using the fact that $J$ is an equivalence of categories, there exists $C' \in \C'$ such that $C \simeq J(C')$. We have:
\begin{align*}
\alpha_{K(C)} \circ \eta_C & = \alpha_{K \circ J(C')} \circ \eta_{J(C')} = \alpha_{K \circ J(C')} \circ G \circ K(\nu_{J(C')}) \circ \gamma_{P(J(C'))} \circ F(\nu^{-1}_{J(C')}).
\end{align*}
On the other hand, we have that $\nu_{J(C')} \circ J(\mu_{C'}) = {\rm id}_{J(C')}$, and so $\nu_{J(C')} = J(\mu^{-1}_{C'})$. Thus, we have:
\begin{align*}
\alpha_{K(C)} \circ \eta_C & = \alpha_{K \circ J(C')} \circ (G \circ K \circ J)(\mu^{-1}_{C'}) \circ \gamma_{P(J(C'))} \circ (F \circ J)(\mu_{C'}) \\
& = \alpha_{K \circ J(C')} \circ \gamma_{C'} \circ (F \circ J)(\mu^{-1}_{C'}) \circ (F \circ J)(\mu_{C'}) & \mbox{(since $\gamma$ is natural)} \\
& = \alpha_{K \circ J(C')} \circ \gamma_{C'} = \eta'_{J(C')} \\
& = \eta'_C.
\end{align*}
We conclude the proof showing the uniqueness of $\alpha$ is the previous expression. Suppose there exists $\alpha' \colon G \Rightarrow G'$ such that $\alpha'_K \circ \eta = \eta'$. Then, after composing with $J$, we obtain  $\alpha'_{K \circ J} \circ \eta_J = \eta'_J$. Using the same arguments as in the previous chain of equalities, we can conclude that $\alpha'_{K \circ J} \circ \gamma = \alpha_{K \circ J} \circ \gamma$, and so $\alpha' = \alpha$. 
\end{myproof}

\begin{corollary}\label{theo:MacLaneExtension}
Suppose that $\C$ is a skeletally small category and that $\E$ is a cocomplete category. Then, for every pair of functors $F \colon \C \too \E$ and $K \colon \C \too \D$, the Kan extension of $F$ along $K$ exists. 
\end{corollary}

\begin{myproof}
Let $\C'$ be a small and dense subcategory of $\C$. Denote the corresponding inclusion by $J \colon \C' \hookrightarrow \C$. On the one hand, since $\C'$ is a full subcategory of $\C$, we have that the functor $J$ is full and faithful. On the other hand, since $\C'$ is dense, we have that $J$ is a \textbf{dense} or \textbf{essentially surjective} functor. Hence, $J$ turns out to be an equivalence of categories. 

Now, since $\C'$ is small, we have that $\Lan_{K \circ J}(F \circ J)$ exists by \emph{Theorem~\ref{theo:MacLane}}. Finally, $\Lan_K(F)$ exists by \emph{Proposition~\ref{prop:MacLaneSke}}.
\end{myproof}

From now on, given a skeletally small preadditive category $\C$, we fix a small and dense subcategory $\C' \subseteq \C$. Due to the previous result, given two functors $F \colon \C \too \Ab$ and $G \colon \C\op \too \Ab$, we can compute $F \otimes_{\C\op} G$ as the tensor product $F \otimes_{(\C')\op} G := \int F \otimes_{\mathbb{Z}} G$, which exists by \emph{Corollary~\ref{coro:existencia_producto}}. Then, in what follows $\C$ will be regarded as a small category. We have the following properties of $F \otimes_{\C\op} G$.

\begin{theorem}
Let $\C$ be a skeletally small preadditive category. For every $F \in \Mod(\C\op)$ and $G \in \Mod(\C)$, the functors
\[
F \displaystyle\operatorname*{\otimes}_{\C\op} - : \Mod(\C) \too \Ab \mbox{ \ and \ } - \displaystyle\operatorname*{\otimes}_{\C\op} G : \Mod(\C\op) \too \Ab
\]
are right exact and preserve arbitrary coproducts in $\Mod(\C)$ and $\Mod(\C\op)$, respectively.
\end{theorem}

\begin{myproof}
We only prove the assertions about $F \otimes_{\C\op} - : \Mod(\C) \too \Ab$, since the proof corresponding to $- \otimes_{\C\op} G : \Mod(\C\op) \too \Ab$ follows by a similar reasoning. 

We first show that $F \otimes_{\C\op} -$ preserves arbitrary coproducts in $\Mod(\C)$. Suppose we have a coproduct $G := \coprod_{i \in I} G_i$ in $\Mod(\C)$. Then, we have:  

\begin{align*}
F \displaystyle\operatorname*{\otimes}_{\C\op} G & = \displaystyle\operatorname*{\int}^{C \in \C\op} F(C) \displaystyle\operatorname*{\otimes}_{\mathbb{Z}} \left( \coprod_{i \in I} G_i(C) \right) \\
& = \displaystyle\operatorname*{\int}^{C \in \C\op} \coprod_{i \in I} \left( F(C) \displaystyle\operatorname*{\otimes}_{\mathbb{Z}} G_i(C) \right) & \mbox{(since $- \displaystyle\operatorname*{\otimes}_{\mathbb{Z}} -$ preserves coproducts)} \\
& = \coprod_{i \in I} F \displaystyle\operatorname*{\otimes}_{\C\op} G_i. & \mbox{(by \emph{Proposition~\ref{prop:cofinal_preservacion}})}
\end{align*}
Hence, $F \otimes_{\C\op} -$ preserves coproducts. 

Now we show that $F \otimes_{\C\op} -$ is right exact. So suppose we are given a short exact sequence in $\Mod(\C)$, say:
\[
0 \Rightarrow G_1 \xRightarrow{\alpha} G_2 \xRightarrow{\beta} G_3 \Rightarrow 0.
\]
Then, for each $C \in \C$, we have the short exact sequence in $\Ab$:
\[
0 \to G_1(C) \xrightarrow{\alpha_C} G_2(C) \xrightarrow{\beta_C} G_3(C) \to 0.
\]
Now, given $C' \in \C$, we know that the functor $F(C') \otimes_{\mathbb{Z}} - : \Ab \to \Ab$ is right exact, which implies that the following sequence is exact in $\Ab$:
\[
F(C') \displaystyle\operatorname*{\otimes}_{\mathbb{Z}} G_1(C) \xrightarrow{F(C') \displaystyle\operatorname*{\otimes}_{\mathbb{Z}} \alpha_C} F(C') \displaystyle\operatorname*{\otimes}_{\mathbb{Z}} G_2(C) \xrightarrow{F(C') \displaystyle\operatorname*{\otimes}_{\mathbb{Z}} \beta_C} F(C') \displaystyle\operatorname*{\otimes}_{\mathbb{Z}} G_3(C) \to 0.
\]
Since $\Ab$ is a Grothendieck category, the direct limit of a direct system of short exact sequences in $\Ab$ (that is, the colimit over a directed set) is also a short exact sequence in $\Ab$. On the other hand, a coproduct of a set of short exact sequences is a particular case of direct limit of its subsets formed by finite coproducts. Thus, we have the following exact sequence in $\Ab$:
\[
\coprod_{f \in \mathsf{Mor}(\C)} F(C') \displaystyle\operatorname*{\otimes}_{\mathbb{Z}} G_1(C) \to \coprod_{f \in \mathsf{Mor}(\C)} F(C') \displaystyle\operatorname*{\otimes}_{\mathbb{Z}} G_2(C) \to \coprod_{f \in \mathsf{Mor}(\C)} F(C') \displaystyle\operatorname*{\otimes}_{\mathbb{Z}} G_3(C) \to 0,
\]
where $f \in \Hom_{\C}(C,C')$. In a similar way, we have the following exact sequence in $\Ab$:
\[
\coprod_{C \in \C} F(C) \displaystyle\operatorname*{\otimes}_{\mathbb{Z}} G_1(C) \to \coprod_{C \in \C} F(C) \displaystyle\operatorname*{\otimes}_{\mathbb{Z}} G_2(C) \to \coprod_{C \in \C} F(C) \displaystyle\operatorname*{\otimes}_{\mathbb{Z}} G_3(C) \to 0.
\]
Thus, we obtain the following commutative diagram in $\Ab$:

\[
\begin{tikzpicture}[description/.style={fill=white,inner sep=2pt}]
\matrix (m) [matrix of math nodes, row sep=2em, column sep=3em, text height=1.75ex, text depth=0.85ex]
{ 
\displaystyle\operatorname*{\coprod}_{f \in \mathsf{Mod}(\C)} F(C') \displaystyle\operatorname*{\otimes}_{\mathbb{Z}} G_1(C) & \displaystyle\operatorname*{\coprod}_{C \in \C} F(C) \displaystyle\operatorname*{\otimes}_{\mathbb{Z}} G_1(C) & F \displaystyle\operatorname*{\otimes}_{\C\op} G_1 & 0 \\
{} & {} & {} & {} \\
\displaystyle\operatorname*{\coprod}_{f \in \mathsf{Mod}(\C)} F(C') \displaystyle\operatorname*{\otimes}_{\mathbb{Z}} G_2(C) & \displaystyle\operatorname*{\coprod}_{C \in \C} F(C) \displaystyle\operatorname*{\otimes}_{\mathbb{Z}} G_2(C) & F \displaystyle\operatorname*{\otimes}_{\C\op} G_2 & 0 \\
{} & {} & {} & {} \\
\displaystyle\operatorname*{\coprod}_{f \in \mathsf{Mod}(\C)} F(C') \displaystyle\operatorname*{\otimes}_{\mathbb{Z}} G_3(C) & \displaystyle\operatorname*{\coprod}_{C \in \C} F(C) \displaystyle\operatorname*{\otimes}_{\mathbb{Z}} G_3(C) & F \displaystyle\operatorname*{\otimes}_{\C\op} G_3 & 0 \\
0 & 0 & 0 \\
};
\path[->]
(m-1-1) edge node[above] {\footnotesize$p_1$} (m-1-2) (m-1-2) edge (m-1-3) (m-1-3) edge (m-1-4)
(m-3-1) edge node[above] {\footnotesize$p_2$} (m-3-2) (m-3-2) edge (m-3-3) (m-3-3) edge (m-3-4)
(m-5-1) edge node[above] {\footnotesize$p_3$} (m-5-2) (m-5-2) edge (m-5-3) (m-5-3) edge (m-5-4)
(m-1-1) edge node[left] {\footnotesize$\displaystyle\operatorname*{\coprod}_{f \in \Hom_{\C}(C',C)} F(C') \displaystyle\operatorname*{\otimes}_{\mathbb{Z}} \alpha_C$} (m-3-1)
(m-3-1) edge node[left] {\footnotesize$\displaystyle\operatorname*{\coprod}_{f \in \Hom_{\C}(C',C)} F(C') \displaystyle\operatorname*{\otimes}_{\mathbb{Z}} \beta_C$} (m-5-1)
(m-1-2) edge node[left] {\footnotesize$\displaystyle\operatorname*{\coprod}_{C \in \C} F(C) \displaystyle\operatorname*{\otimes}_{\mathbb{Z}} \alpha_C$} (m-3-2)
(m-3-2) edge node[left] {\footnotesize$\displaystyle\operatorname*{\coprod}_{C \in \C} F(C) \displaystyle\operatorname*{\otimes}_{\mathbb{Z}} \beta_C$} (m-5-2)
(m-1-3) edge node[left] {\footnotesize$F \displaystyle\operatorname*{\otimes}_{\C\op} \alpha$} (m-3-3)
(m-3-3) edge node[left] {\footnotesize$F \displaystyle\operatorname*{\otimes}_{\C\op} \beta$} (m-5-3)
(m-5-1) edge (m-6-1) 
(m-5-2) edge (m-6-2)
(m-5-3) edge (m-6-3)
;
\end{tikzpicture}
\]
where: 
\begin{align*}
p_1 & := \coprod_{f \in \mathsf{Mor}(\C)} F(D') \otimes G_1(f) - \coprod_{f \in \mathsf{Mor}(\D)} F(f) \otimes G_1(D), \\ 
p_2 & := \coprod_{f \in \mathsf{Mor}(\C)} F(D') \otimes G_2(f) - \coprod_{f \in \mathsf{Mor}(\D)} F(f) \otimes G_2(D), \\ 
p_3 & := \coprod_{f \in \mathsf{Mor}(\C)} F(D') \otimes G_3(f) - \coprod_{f \in \mathsf{Mor}(\D)} F(f) \otimes G_3(D).
\end{align*} 
Note that the rows and the first two columns (from left to right) in the previous diagram, are exact sequences. Since the bottom-right square commutes, it is easy to check that the morphism $F \otimes_{\C\op} \beta$ is epic. Therefore, the right-hand column is an exact sequence in $\Ab$, that is, the functor $F \otimes_{\C\op} -$ is right exact. 
\end{myproof}

We conclude these notes proving, in two parts, that $F \otimes_{\C\op} G = F \otimes G$, that is, we first show that the previous equality holds in the case $G$ is free, and then by using free resolutions for $G$ arbitrary.

Let $G = \coprod_{i \in I} \Hom_{\C}(-,C_i)$ be a free $\C$-module. For each $G_i = \Hom_{\C}(-,C_i)$, we know that: 
\begin{align*}
& F \displaystyle\operatorname*{\otimes}_{\C^{\rm op}} G_i = F(C_i), & \mbox{(by the dual of \eqref{eqn:derecha})} \\
& F \otimes G_i = F(C_i). & \mbox{(by \emph{Theorem~\ref{theo:producto_tensorial}})}
\end{align*}
Since $F \otimes_{\C\op} -$ and $F \otimes -$ preserve coproducts, we have:
\begin{align*}
& F \displaystyle\operatorname*{\otimes}_{\C^{\rm op}} G = \coprod_{i \in I} F \displaystyle\operatorname*{\otimes}_{\C^{\rm op}} G_i = \coprod_{i \in I} F(C_i) = \coprod_{i \in I} F \otimes G_i = F \otimes G.
\end{align*}
Thus, $F \otimes_{\C\op} G \simeq F \otimes G$, for every free $\C$-module $G$. 

Now let $G \in \Mod(\C)$. Consider a partial free resolution of $G$, say:
\[
\coprod_{i \in I} \Hom_{\C}(-,C_i) \Rightarrow \coprod_{j \in J} \Hom_{\C}(-,C_j) \Rightarrow G \Rightarrow 0.
\]
Applying $F \otimes_{\C\op} -$ and $F \otimes -$, which are right exact and preserve arbitrary coproducts, we obtain the following commutative diagram in $\Ab$ with exact rows:
\[
\begin{tikzpicture}[description/.style={fill=white,inner sep=2pt}]
\matrix (m) [matrix of math nodes, row sep=3em, column sep=3em, text height=1.25ex, text depth=0.25ex]
{ 
F \displaystyle\operatorname*{\otimes}_{\C\op} \displaystyle\operatorname*{\coprod}_{i \in I} \Hom_{\C}(-,C_i) & F \displaystyle\operatorname*{\otimes}_{\C\op} \displaystyle\operatorname*{\coprod}_{j \in J} \Hom_{\C}(-,C_j) & F \displaystyle\operatorname*{\otimes}_{\C\op} G & 0 \\
F \otimes \displaystyle\operatorname*{\coprod}_{i \in I} \Hom_{\C}(-,C_i) & F \otimes \displaystyle\operatorname*{\coprod}_{j \in J} \Hom_{\C}(-,C_j) & F \otimes G & 0 \\
};
\path[->]
(m-1-1) edge (m-1-2) (m-1-2) edge (m-1-3) (m-1-3) edge (m-1-4)
(m-2-1) edge (m-2-2) (m-2-2) edge (m-2-3) (m-2-3) edge (m-2-4)
;
\path[-,font=\scriptsize]
(m-1-1) edge [double, thick, double distance=2pt] (m-2-1)
(m-1-2) edge [double, thick, double distance=2pt] (m-2-2)
;
\path[dotted,->]
(m-1-3) edge node[left] {\footnotesize$\exists!$} node[right] {\footnotesize$h$} (m-2-3)
;
\end{tikzpicture}
\]
where the morphism $h$, induced by the universal property of coproducts, turns out to be an isomorphism. Therefore,
\[
F \displaystyle\operatorname*{\otimes}_{\C\op} G \simeq F \otimes G,
\]
for every $F \in \Mod(\C\op)$ and $G \in \Mod(\C)$. The fact that the previous isomorphism in natural is easy to check, and it is left to the reader.


\section*{Acknowldgements}

These notes are the result of a mini-seminar on Kan extensions taught at the Instituto de Matem\'aticas in the Universidad Nacional Aut\'onoma de M\'exico (IMATE-UNAM), between May and June of 2016. The author wants to thank the audience for their interest and participation. Among them, to Octavio Mendoza, Violeta Garc\'ia, Yadira Valdivieso, V\'ictor Becerril, Mindy Huerta, Luis Mart\'inez and Alejandro Argud\'in, from the IMATE-UNAM; and to Edith Corina S\'aenz, Valente Santiago, and Andr\'es Barei, from the UNAM Faculty of Sciences; and to Mart\'in Ortiz-Morales, from the Universidad Aut\'onoma del Estado de M\'exico. Special thanks to Violeta Garc\'ia and Octavio Mendoza, for answering some of my doubts on the topics presented here, and for recommending me several useful references. Special thanks go also to Omar Antol\'in, from the University of British Columbia, who made several remarks and corrections regarding the terminology employed in previous versions of these notes, which have improved the quality of the current version.

\newpage


\addcontentsline{toc}{section}{Referencias}
\bibliographystyle{alpha}
\bibliography{BiblioKan}

\end{document}